\documentclass[11pt]{amsart}
\newtheorem{theorem}{Theorem}[section]

\newtheorem{lemma}[theorem]{Lemma}
\newtheorem{proposition}[theorem]{Proposition}

\newtheorem{example}[theorem]{Example}
\theoremstyle{definition}
\newtheorem{definition}[theorem]{Definition}
\theoremstyle{remark}
\newtheorem{remark}[theorem]{Remark}
\numberwithin{equation}{section}

\newcommand{\CC}{\mathcal C}
\newcommand{\B}{\mathbb B}
\newcommand{\C}{\mathbb C}
\newcommand{\R}{\mathbb R}

\begin{document}

\title[Fefferman's mapping theorem on almost complex manifolds]
{Fefferman's mapping theorem on almost complex manifolds}

\author{Bernard Coupet, Herv\'e Gaussier and Alexandre Sukhov}

\address{\begin{tabular}{lllll}
Bernard Coupet & & Herv\'e Gaussier & & Alexandre Sukhov\\
C.M.I. & & C.M.I. & & U.S.T.L. \\
39, rue Joliot-Curie, & & 39, rue Joliot-Curie, & & Cit\'e Scientifique \\
13453 Marseille Cedex 13 & & 13453 Marseille Cedex 13 & &
59655 Villeneuve d'Ascq Cedex\\
FRANCE & & FRANCE & & FRANCE\\
& & & & \\
{\rm coupet@cmi.univ-mrs.fr} & & {\rm gaussier@cmi.univ-mrs.fr} & &
{\rm sukhov@agat.univ-lille1.fr}
\end{tabular}
}


\subjclass[2000]{32H02, 53C15}

\date{\number\year-\number\month-\number\day}

\begin{abstract}
We give a necessary and sufficient condition for the smooth extension of
a diffeomorphism between smooth strictly pseudoconvex domains in four
real dimensional almost complex manifolds (see Theorem~\ref{MainTheorem}).
The proof is mainly based on a reflection principle for
pseudoholomorphic discs, on precise estimates of the Kobayashi-Royden
infinitesimal pseudometric and on the scaling method in almost complex
manifolds.
\end{abstract}

\maketitle

\section{Introduction}

The analysis on almost complex manifolds, first developed by
Newlander-Nirenberg and Nijenhuis-Woolf, appeared crucial in
symplectic and contact geometry with the fundamental work of M.Gromov
\cite{Gr}. Since the literature dedicated to this subject is rapidly
growing we just mention the book \cite{AuLa} and references
therein. The (geometric) analysis on almost complex
manifolds became one of the most powerful tools in symplectic
geometry, making its systematic developpement relevant. The present
paper is a step in this program.

Fefferman's mapping theorem~\cite{Fe} states that a biholomorphism between two
smoothly bounded strictly pseudoconvex domains in $\C^n$ extends as a smooth
diffeomorphism between their closures. This result had a strong impact on the
developpement of complex analysis on domains in $\C^n$ during the
last twenty five years. Our main goal is to prove an analogue of 
this theorem in almost complex manifolds. Complex and symplectic structures
are usually related as follows. Let $(M,\omega)$ and $(M',\omega')$ be
two real manifolds equipped with symplectic forms and let $J$ be an almost
complex structure on $M$ tamed by $\omega$ (so that $\omega(v,Jv)
> 0$ for any non-zero vector $v$). If $\phi: M \rightarrow M'$ is
a symplectomorphism,  the direct image $J':= \phi_*(J) = d\phi \circ J
\circ d\phi^{-1}$ of the structure $J$ is an almost complex structure 
on $M'$ tamed by
$\omega'$ and $\phi$ is a biholomorphism with respect to $J$ and
$J'$. This property enables to construct
topological invariants of symplectic structures employing the complex
geometry. In his survey \cite{Be}, D.Bennequin raised the question of a
symplectic analogue of this theorem.
E.Chirka constructed in \cite{Ch} an example of a symplectomorphism
of the unit ball in $\C^n$ with the usual symplectic structure, having
a wild boundary behaviour. This gives a negative answer to the
question. Our main result shows that Fefferman's theorem remains true
in the category of almost complex manifolds~:

\begin{theorem}
\label{MainTheorem}
Let $D$ and $D'$ be two smooth relatively compact domains in real
four dimensional  manifolds. Assume that
$D$ admits an almost complex structure $J$ smooth on $\bar D$ and
 such that $(D,J)$ is strictly pseudoconvex. Then a smooth
diffeomorphism $f:  D  \rightarrow  D'$ extends to a smooth
diffeomorphism between $\bar D$ and $\bar D'$  if and only
if the direct image $f_*(J)$ of $J$ under
$f$ extends smoothly on $ \bar D'$ and $(D', f_*(J))$ 
is strictly pseudoconvex.
\end{theorem}

One can see the smooth extension of the direct image $f_*(J)$ to the boundary
of $D'$ as the smooth extension, up to the boundary, of a part 
of first order partial derivatives of the components of $f$.
Theorem~\ref{MainTheorem} claims that all the partial derivatives 
necessarily extend smoothly up to the boundary. This statement is a geometric
version of the elliptic regularity and is a criterion applicable (at least
in principle) to any diffeomorphism between four dimensional
real manifolds with boundaries. 

Theorem~\ref{MainTheorem} admits the following formulation, closer to the
claasical one.

\begin{theorem}\label{TH2}
A biholomorphism between two smooth relatively compact strictly pseudoconvex
domains in (real) four dimensional almost complex 
manifolds extends to a smooth diffeomorphism between their closures.
\end{theorem}

Theorem~\ref{TH2} is a generalization of Fefferman's mapping theorem in
(complex) dimension 2. We point out that in the almost complex
category, the real four dimensional manifolds often represent the most
interesting class from the point of view of
symplectic geometry. Our restriction on the dimension comes from our method
of proof; we do not know if it necessary in the general case.

The original proof of C.Fefferman is based on a subtil
investigation of asymptotic behavior of the Bergman kernel in strictly
pseudoconvex domains. Later several different approaches have been
proposed. Similarly to the integrable case (see, for instance,
L.Nirenberg-S.Webster-P.Yang \cite{NWY}, S.Pinchuk-S.Khasanov \cite{PiKh},
B.Coupet \cite{Co} and F.Forstneric \cite{Fo})
our proof is based on boundary estimates of the infinitesimal Kobayashi-
Royden pseudometric and on the smooth reflection principle.
We point out that
in the almost complex case, the reflection
principle for totally real manifolds has been used by H.Hofer
\cite{Ho}, S.Ivashkovich-V.Shevchishin \cite{IvSh}, E.Chirka
\cite{Ch1}.

The paper is organized as follows. 

Section 2 is preliminary and contains
general facts about almost complex manifolds. In Subsection~2.5 
for a given point on a strictly pseudoconvex hypersurface $\Gamma$
in a four dimensional almost complex manifold $(M,J)$, we
construct a coordinates chart in which $J$ is a sufficiently small
diagonal perturbation of the standard structure $J_{st}$ on $\C^2$
and $\Gamma$ is strictly pseudoconvex with respect to $J_{st}$.
Such a representation will be crucially used and explains the
restriction to the four dimensional case in our approach.

In Section~3 we prove that in the hypothesis of Theorem~\ref{MainTheorem}
a biholomorphism $f$ extends as a $1/2$-H\"older map between the closures
of the domains and we study the boundary behaviour of its tangent map.
Similarly to the integrable case, our proof is based on the Hopf lemma
and the estimates of the Kobayashi-Royden infinitesimal pseudometric
obtained in \cite{GaSu}. This result allows to restrict our considerations
to the case where $f$ is a biholomorphism between two strictly pseudoconvex
domains in $\C^2$ with small almost complex deformations of the standard
structure.

Sections~4 and 5 contain another technical ingredient necessary for the proof
of Fefferman's theorem : results on the boundary regularity of
pseudoholomorphic maps near totally real manifolds. In Section~4 we study
the boundary regularity of a pseudoholomorphic disc attached, in the sense
of the cluster set, to a totally real submanifold of an almost complex
manifold. The proof consists of two steps. First we obtain an a priori
bound for the gradient of the disc using uniform estimates of the
Kobayashi-Royden infinitesimal pseudometric in the Grauert tube around a
totally real manifold (in the integrable case a similar construction has
been used in \cite{ChCoSu}).
Then we apply a version of the smooth reflection principle in almost complex
manifolds (this construction is due to E.Chirka~\cite{Ch1}).
In Section~5 we establish the boundary regularity of a pseudoholomorphic
map defined in a wedge with a totally real edge and taking this edge (in the
sense of the cluster set) to a totally real submanifold in an almost complex
manifold. For the proof we fill the wedge by pseudoholomorphic discs
attached to the edge along the upper semi circle (this generalizes
Pinchuk's construction in the integrable case \cite{Pi74})
and we apply the results of Section~4. 

In Section~6 we show how to deduce the proof of Fefferman's theorem
from the results of Section~5. The main idea is to consider the cotangent
lift $\tilde f$ of a biholomorphism $f$. According to
the known results of the differential geometry \cite{YI}, an
almost complex structure $J$ on $M$ admits a canonical almost complex
lift $\tilde J$ to the cotangent bundle of $M$ such that $\tilde f$ is
biholomorphic with respect to $\tilde J$ and the conormal bundle of
a strictly $J$-pseudoconvex hypersurface is totally real with
respect to $\tilde J$. In the integrable case
the holomorphic tangent bundle of a strictly
pseudoconvex hypersurface (that is the projectivization of the conormal
bundle) is frequently used instead of the conormal bundle
(see \cite{We, PiKh, Co}). Since in the almost complex case the
projectivization of the cotangent bundle does not admit a natural
almost complex structure, we need to deal with the conormal bundle,
similarly to the ideas of A.Tumanov \cite{Tu94}.
In the case where $f$ is of class $\mathcal C^1$ up to the boundary of
$D$, its cotangent lift extends continuously to the conormal bundle of
$\partial D$ and takes it to the conormal bundle of $\partial D'$, which
implies Fefferman's theorem in that case in view of the results of
Section~5.

In Section~7 we consider the general situation of Theorem~\ref{MainTheorem}. 
We prove that the cotangent lift of $f$ takes the conormal bundle of
$\partial D$ to the conormal bundle of $\partial D'$ in the sense of
the cluster set. This is sufficient in order to apply the results of
Section~5. Our proof is based on the results of Section~3 and the scaling
method introduced by S.Pinchuk \cite{Pi} in the integrable case which we
develop in our situation.

\section{Preliminaries}

\subsection{Almost complex manifolds.}

Let $(M',J')$ and $(M,J)$ be almost complex manifolds and let $f$ be
a smooth map from $M'$ to $M$. We say that $f$ is {\sl
 $(J',J)$-holomorphic} if $df \circ J' = J \circ df$ on $TM'$. 
We denote by $\mathcal O_{(J',J)}(M',M)$ the set of
$(J',J)$-holomorphic maps from $M'$ to $M$. Let $\Delta$ be the unit
disc in $\C$ and $J_{st}$ be the standard integrable structure on $\C^n$
for every $n$. If $(M',J')=(\Delta,J_{st})$,
we denote by $\mathcal O_J(\Delta,M)$ the set
$\mathcal O_{(J_{st},J)}(\Delta,M)$ of {\sl $J$-holomorphic discs} in $M$. 

The following Lemma shows that every almost complex manifold
$(M,J)$ can be viewed locally as the unit ball in
$\mathbb C^n$ equipped with a small almost complex
deformation of $J_{st}$. This will be used frequently in the sequel.
\begin{lemma}
\label{suplem1}
Let $(M,J)$ be an almost complex manifold. Then for every point $p \in
M$ and every $\lambda_0 > 0$ there exist a neighborhood $U$ of $p$ and a
coordinate diffeomorphism $z: U \rightarrow \mathbb B$ such that
$z(p) = 0$, $dz(p) \circ J(p) \circ dz^{-1}(0) = J_{st}$  and the
direct image $\hat J = z_*(J)$ satisfies $\vert\vert \hat J - J_{st}
\vert\vert_{\CC^2(\bar {\mathbb B})} \leq \lambda_0$.
\end{lemma}
\proof There exists a diffeomorphism $z$ from a neighborhood $U'$ of
$p \in M$ onto $\mathbb B$ satisfying $z(p) = 0$ and $dz(p) \circ J(p)
\circ dz^{-1}(0) = J_{st}$. For $\lambda > 0$ consider the dilation
$d_{\lambda}: t \mapsto \lambda^{-1}t$ in $\C^n$ and the composition
$z_{\lambda} = d_{\lambda} \circ z$. Then $\lim_{\lambda \rightarrow
0} \vert\vert (z_{\lambda})_{*}(J) - J_{st} \vert\vert_{\CC^2(\bar
{\mathbb B})} = 0$. Setting $U = z^{-1}_{\lambda}(\mathbb B)$ for
$\lambda > 0$ small enough, we obtain the desired statement. \qed

\subsection{ $\partial_J$ and $\bar{\partial}_J$ operators}

Let $(M,J)$ be an almost complex manifold. We denote by $TM$ the real 
tangent bundle of $M$ and by $T_\C M$ its complexification. Recall
that $T_\C M = T^{(1,0)}M \oplus T^{(0,1)}M$ where
$T^{(1,0)}M:=\{ X \in T_\C M : JX=iX\} = \{\zeta -iJ \zeta, \zeta \in
TM\},$ 
and $T^{(0,1)}M:=\{ X \in T_\C M : JX=-iX\} = \{\zeta +
iJ \zeta, \zeta \in TM\}$.
 Let $T^*M$ denote the cotangent bundle of  $M$.
Identifying $\C \otimes T^*M$ with
$T_\C^*M:=Hom(T_\C M,\C)$ we define the set of complex
forms of type $(1,0)$ on $M$ by~:
$
T^*_{(1,0)}M=\{w \in T_\C^* M : w(X) = 0, \forall X \in T^{(0,1)}M\}
$
and the set of complex forms of type $(0,1)$ on $M$ by~:
$
T^*_{(0,1)}M=\{w \in T_\C^* M : w(X) = 0, \forall X \in T^{(1,0)}M\}
$.
Then $T_\C^*M=T^*_{(1,0)}M \oplus T^*_{(0,1)}M$.   
This allows to define the operators $\partial_J$ and
$\bar{\partial}_J$ on the space of smooth functions defined on
$M$~: given a complex smooth function $u$ on $M$, we set $\partial_J u =
du_{(1,0)} \in T^*_{(1,0)}M$ and $\bar{\partial}_Ju = du_{(0,1)}
\in T^*_{(0,1)}M$. As usual,
differential forms of any bidegree $(p,q)$ on $(M,J)$ are defined
by means of the exterior product.

\subsection{Real submanifolds in an almost complex manifold}

Let $\Gamma$ be a real smooth submanifold in $M$ and let $p \in
\Gamma$. We denote by $H^J(\Gamma)$ the $J$-holomorphic tangent bundle
$T\Gamma \cap JT\Gamma$. 

\begin{definition}
The real submanifold $\Gamma$ is totally real if $H^J(\Gamma)=\{0\}$.
\end{definition}

We note that if $\Gamma$ is a real
hypersurface in $M$ defined by $\Gamma=\{r=0\}$ and $p \in \Gamma$
then by definition $H_p^J(\Gamma) = \{v \in T_pM : dr(p)(v) =
dr(p)(J(p)v) = 0\}$. 

We recall the notion of the Levi form of a hypersurface~:

\begin{definition}\label{DEF}
Let $\Gamma=\{r=0\}$ be a smooth real hypersurface in $M$ 
($r$ is any smooth defining function of $\Gamma$) and let $p \in \Gamma$. 

$(i)$ The {\sl Levi form} of $\Gamma$ at $p$ is the map defined on
$H^J_p(\Gamma)$  by ${\mathcal L}_{\Gamma}^J(X_p) = J^\star dr[X,JX]_p$,
where the vector field $X$ is any section of the $J$-holomorphic tangent
bundle  $H^J \Gamma$ such that $X(p) = X_p$.

$(ii)$ A real smooth hypersurface $\Gamma=\{r=0\}$ in $M$ is 
{\sl strictly $J$-pseudoconvex} if its Levi form ${\mathcal L}_{\Gamma}^J$ 
is positive definite on $H^J(\Gamma)$.
\end{definition}
  
\begin{remark} 
$(i)$ the ``strict $J$-pseudoconvexity'' condition does not depend on
the choice of a smooth defining function of $\Gamma$. Indeed if $\rho$
is an other smooth defining function for $\Gamma$ in a neighborhood of
$p \in \Gamma$ then there exists a positive smooth function $\lambda$
defined in a neighborhood of $p$ such that $\rho=\lambda r$.  In
particular $(J^\star dr)(p) = \lambda(p)(J^\star d\rho)(p)$.

$(ii)$ since the map $(r,J) \mapsto J^\star dr$ is smooth the ``strict
$J$-pseudoconvexity'' is stable under small perturbations of both the
hypersurface and the almost complex structure.
\end{remark}

Let $X \in TM$. It follows from the identity
$d(J^\star dr)(X,JX)=X(<J^\star dr,JX>) - JX(<J^\star dr,X>) - 
(J^\star dr)[X,JX]$
that 
$
(J^\star dr)[X,JX] = -d(J^\star dr)(X,JX)
$ 
for every $X \in H^J\Gamma$, since $<dr,JX>=<dr,JX> = 0$ in that case. 
Hence we set

\begin{definition} If $r$ is a $\CC^2$ function on $M$ then the Levi 
form of $r$ is defined on $TM$ by ${\mathcal L}^J(r)(X):=
-d(J^\star dr)(X,JX)$.
\end{definition}
Let $p \in M$ and $v \in T_pM$. We will denote by ${\mathcal L}^J(r)(p)(v)$ 
the quantity ${\mathcal L}^J(r)(X)(p)$ where $X$ is any section of $TM$ 
such that $X(p) = v$. Obviously, the Levi form ${\mathcal L}^J(r)$ is
determined by the  form $-d(J^\star dr)_{(1,1)}$ (where the  $(1,1)$
part of $-d(J^\star dr)(X,JX)$ is taken with respect to $J$).

\subsection{Kobayashi-Royden infinitesimal pseudometric}
Let $(M,J)$ be an almost complex
manifold. In what follows we use the notation $\zeta=x+iy \in \mathbb C$.
According to \cite{NiWo}, for every $p \in M$ there is a
neighborhood $\mathcal V$ of $0$ in $T_pM$ such that for every $v \in
\mathcal V$ there exists $f \in \mathcal O_J(\Delta,M)$ satisfying
$f(0) = p,$ $df(0) (\partial / \partial x) = v$. This allows to define
the Kobayashi-Royden infinitesimal pseudometric $K_{(M,J)}$.
\begin{definition}\label{dd}
For $p \in M$ and $v \in T_pM$, $K_{(M,J)}(p,v)$ is the infimum of the
set of positive $\alpha$ such that there exists a $J$-holomorphic disc
$f:\Delta \rightarrow M$ satisfying $f(0) = p$ and $df(0)(\partial
/\partial x) = v/\alpha$.
\end{definition}
Since for every $f \in \mathcal O_{(J',J)}(M',M)$ and
every $\varphi \in \mathcal O_J(\Delta,M')$ the composition $f \circ
\varphi$ is in $\mathcal O_J(\Delta,M)$ we have~:
\begin{proposition}\label{ppp}
Let $f:(M',J') \rightarrow (M,J)$ be a $(J',J)$-holomorphic map. 
Then $K_{(M,J)}(f(p'),df(p')(v')) \leq K_{(M',J')}(p',v')$ for every
$p'\in M', \ v' \in T_{p'}M'$.
\end{proposition}

\subsection{Plurisubharmonic functions.} We first recall the following 
definition~:
\begin{definition}\label{d6}
An upper semicontinuous function $u$ on $(M,J)$ is called 
{\sl $J$-plurisubharmonic} on $M$ if the composition $u \circ f$ 
is subharmonic on $\Delta$ for every $f \in \mathcal O_J(\Delta,M)$.
\end{definition} 
If $M$ is a domain in $\C^n$ and $J=J_{st}$ then a 
$J_{st}$-plurisubharmonic function is a plurisubharmonic function 
in the usual sense. 

The next proposition gives a characterization of
$J$-plurisubharmonic functions (see \cite{de99,ha02})~:
\begin{proposition}\label{PROP}
Let $u$ be a $\CC^2$ real valued function on $M$. Then $u$ is
$J$-plurisubharmonic on $M$ if and only if $\mathcal L^J(u)(X) \geq 0$
for every $X \in TM$.
\end{proposition}
Proposition~\ref{PROP} leads to the definition~:
\begin{definition}
A $\CC^2$ real valued function $u$ on $M$ is {\sl strictly 
$J$-plurisubharmonic} on $M$ if  $\mathcal L^J(u)$
is positive definite on $TM$.
\end{definition}

We have the following example of a
$J$-plurisubharmonic function on an almost complex manifold $(M,J)$~:
\begin{example}\label{example}
For every point $p\in (M,J)$ there exists a neighborhood $U$ of $p$
 and a diffeomorphism $z:U \rightarrow \mathbb
B$ centered at $p$ (ie $z(p) =0$) such that the function $|z|^2$ is
$J$-plurisubharmonic on $U$. 
\end{example}
\proof Let $p \in M$, $U_0$ be a neighborhood of $p$ and $z: U_0
\rightarrow \mathbb B$ be local complex coordinates centered at $p$,
such that $dz \circ J(p) \circ dz^{-1} = J_{st}$ on $\mathbb B$. Consider
the function $u(q) = |z(q)|^2$ on $U_0$. For every $w,v \in \C^n$ we
have $\mathcal L^{J_{st}}(u\circ z^{-1})(w)(v) = \|v\|^2$. Let
$B(0,1/2)$ be the ball centered at the origin with radius $1/2$ and
let $\mathcal E$ be the space of smooth almost complex structures
defined in a neighborhood of $\overline{B(0,1/2)}$. Since the function
$(J',w) \mapsto \mathcal L^{J'}(u \circ z^{-1})(w)$ is continuous on
$\mathcal E \times B(0,1/2)$, there exist a neighborhood $V$ of the
origin and positive constants $\lambda_0$ and $c$ such that $\mathcal
L^{J'}(u\circ z^{-1})(q)(v)\geq c \|v\|^2$ for every $q \in V$ and
for every almost complex structure $J'$ satisfying
$\|J'-J_{st}\|_{\CC^2(\bar{V})} \leq \lambda_0$. Let $U_1$ be a
neighborhood of $p$ such that $\|z_*(J) -
J_{st}\|_{\CC^2(\overline{z(U_1)})}\leq \lambda_0$ and let $0<r<1$ be
such that $B(0,r) \subset V$ and $U:=z^{-1}(B(0,r)) \subset U_1$. Then
we have the following estimate for every $q \in U$ and $v \in T_qM$~:
$\mathcal L^J(u)(q)(v) \geq c \|v\|^2$. Then $r^{-1}z$ is the desired
diffeomorphism. \qed

\vskip 0,1cm     
We also have the following 

\begin{lemma}
A function  $u$ of class $\CC^2$ in a neighborhood
of a point $p$ of $(M,J)$  is strictly $J$-plurisubharmonic
if and only if there exists a neighborhood $U$ of $p$  with local
complex coordinates $z:U \rightarrow \mathbb B$ centered at $p$, such
that the function $u - c|z|^2$ is $J$-plurisubharmonic on $U$ for some
constant $c > 0$.
\end{lemma}

\subsection{Local description of strictly pseudoconvex domains.} If
$\Gamma$ is a germ of a real hypersurface in $\C^n$ strictly
pseudoconvex with respect to $J_{st}$, then $\Gamma$ remains strictly
pseudoconvex for any almost complex structure $J$ sufficiently close
to $J_{st}$  in the $\CC^2$-norm. 
Conversely a strictly
pseudoconvex hypersurface in an almost complex manifold of real dimension
four can be represented, in suitable local coordinates, as a
strictly $J_{st}$-pseudoconvex hypersurface equipped with a small deformation
of the standard structure. Indeed, according to \cite{Si} Corollary~3.1.2,
there exist a neighborhood $U$ of $q$ in $M$ and complex coordinates
$z=(z^1,z^2) : U \rightarrow  B \subset \C^2$, $z(q) =
0$ such that $z_*(J)(0) = J_{st}$ and moreover, a map $f: \Delta
\rightarrow  B$ is $J':= z_*(J)$-holomorphic if it satisfies the
equations 

\begin{eqnarray}
\label{Jhol}
\frac{\partial f^j}{\partial \bar \zeta} =
A_j(f^1,f^2)\overline{\left ( \frac{\partial f^j}
{\partial \zeta}\right ) }, j=1,2
\end{eqnarray} 
where $A_j(z) =  O(\vert
z \vert)$, $j=1,2$.

In order to obtain such coordinates, one can consider two transversal
foliations of the ball $\mathbb B$ by $J'$-holomorphic curves
(see~\cite{NiWo})and then take these curves into the lines $z^j = const$
by a local diffeomorphism. The direct image of the almost complex structure
$J$ under such a diffeomorphism has a diagonal matrix $ J'(z^1,z^2) =
(a_{jk}(z))_{jk}$ with $a_{12}=a_{21}=0$ and $a_{jj}=i+\alpha_{jj}$
where $\alpha_{jj}(z)=\mathcal O(|z|)$ for $j=1,2$.
We point out that the lines $z^j = const$ are
$J$-holomorphic after a suitable parametrization (which, in general,
is not linear). 

In what follows we omit the prime and denote this structure again by
$J$. We may assume that the complex tangent space $T_0(\partial D)
\cap J(0) T_0(\partial D) = T_0(\partial D) \cap i T_0(\partial D)$ is
given by $\{ z^2 = 0 \}$.
In particular, we have the following expansion for the defining
function $\rho$ of $D$ on $U$~:
$\rho(z,\bar{z}) = 2 Re(z^2) + 2Re K(z) + H(z) + \mathcal O(\vert z
\vert^3)$, where
$K(z)  = \sum k_{\nu\mu} z^{\nu}{z}^{\mu}$, $k_{\nu\mu} =
k_{\mu\nu}$ and 
$H(z) = \sum h_{\nu\mu} z^{\nu}\bar z^{\mu}$, $h_{\nu\mu} =
\bar h_{\mu\nu}$.

\begin{lemma}
\label{PP}
The domain $D$  is strictly $J_{st}$-pseudoconvex near the origin.
\end{lemma}

\noindent{\it Proof of Lemma~\ref{PP}.} Consider a complex vector
 $v=(v_1,0)$ tangent to $\partial D$ at the origin.  Let $f:\Delta
 \rightarrow \C^2$ be a $J$-holomorphic disc centered at the
 origin and tangent to $v$: $f(\zeta) = v\zeta + \mathcal O(\vert
 \zeta \vert^2)$.  Since $A_2 = \mathcal O(\vert z \vert)$, it follows
 from the $J$-holomorphy equation (\ref{Jhol}) that
 $(f^2)_{\zeta\bar\zeta}(0) = 0$. This implies that $(\rho \circ
 f)_{\zeta\bar\zeta}(0) = H(v).$ Thus, the Levi form with respect to
 $J$ coincides with the Levi form with respect to $J_{st}$ on the
 complex tangent space of $\partial D$ at the origin. This proves
Lemma~\ref{PP}. \qed

\vskip 0,1cm
Consider the non-isotropic dilations $\Lambda_{\delta}: (z^1,z^2) \mapsto
(\delta^{-1/2}z^1,\delta^{-1}z^2) = (w^1,w^2)$ with $\delta > 0$. 
If $J$ has the above
diagonal form in the coordinates $(z^1,z^2)$ in $\C^2$, then
its direct image  $J_{\delta}= (\Lambda_{\delta})_*(J)$ has the form
$J_{\delta}(w^1,w^2) =(a_{jk}(\delta^{1/2}w^1,\delta w^2))_{jk}$
and so $J_{\delta}$ tends to $J_{st}$ in the $\mathcal C^2$ norm as $\delta
\rightarrow 0$. On the other hand, $\partial D$ is, in the $w$ coordinates,
the zero set of the function 
$\rho_{\delta}= \delta^{-1}(\rho \circ \Lambda_{\delta}^{-1})$.
As $\delta \rightarrow 0$, the function $\rho_{\delta}$ tends to 
the function $2 Re w^2 + 2 Re K(w^1,0) + H(w^1,0)$ which defines a
strictly $J_{st}$-pseudoconvex domain by Lemma~\ref{PP} and proves the claim.

This also proves that if $\rho$
is a local defining function of a strictly $J$-pseudoconvex domain, then
$\tilde{\rho}:=\rho + C \rho^2$
is a strictly $J$-plurisubharmonic function, quite similarly to the standard
case.

In conclusion we point out that extending $\tilde \rho$ by a suitable
negative constant, we obtain that if $D$ is a strictly
$J$-pseudoconvex domain in an almost complex
manifold, then there exists a neighborhood $U$ of $\bar{D}$ and a
function $\rho$, $J$-plurisubharmonic on $U$ and strictly
$J$-plurisubharmonic in a neighborhood of $\partial D$, such that
$D=\{ \rho <0\}$.
 
\section{ Boundary continuity and localization of biholomorphisms}
 
In this section we give some preliminary technical results necessary
for the proof of Theorem~\ref{MainTheorem}.

\subsection{Hopf lemma and the  boundary distance preserving property}

In what follows we need an analog of the Hopf lemma 
for almost complex manifolds. It can be proved quite similarly 
to the standard one. 

\begin{lemma}\label{hopf}
(Hopf lemma) Let $G$ be a relatively compact domain with a $\mathcal C^2$
boundary on an almost complex manifold $(M,J)$. Then for any negative
$J$-psh function $u$ on $D$ there exists a constant $C > 0$ such that
$\vert u(p) \vert \geq C dist(p,\partial G)$ for any $p \in G$ ($dist$
is taken with respect to a Riemannian metric on $M$).
\end{lemma}
\noindent{\it Proof of Lemma~\ref{hopf}}.
{\it Step 1.} We have the following precise version on the unit
disc: let $u$ be a subharmonic function on $\Delta$, $K$ be a fixed
compact on $\Delta$. Suppose that $u < 0$ on $\Delta$ and $u \vert K
\leq -L$ where $ L > 0$ is constant. Then there exists $C(K,L) > 0$
(independent of $u$) such that $\vert u(p) \vert \geq C
dist(p,\partial \Delta)$ (see \cite{Ra}). 

{\it Step 2.} Let $G$ be a domain in $\C$ with $\mathcal C^2$-boundary.
Then there exists an $r > 0$ (depending on the curvature of the boundary)
such that for any boundary point $q \in \partial G$ the ball $B_{q,r}$
of radius $r$ centered on
the interior normal to $\partial G$ at $q$, such that  $q \in
\partial B_{q,r}$, is
contained in $G$. Applying Step 1 to the restriction of $u$ on every
such a ball (when $q$ runs over $\partial G$) we obtain the Hopf lemma
for a domain with $\CC^2$ boundary: 
let $u$ be a subharmonic function on $G$, $K$ be a fixed
compact on $G$. Suppose that $u < 0$ on $G$ and $u \vert K
\leq -L$ where $ L > 0$ is constant. Denote by $k$ the curvature of
$\partial G$. Then there exists $C(K,L,k) > 0$
(independent of $u$) such that $\vert u(p) \vert \geq C
dist(p,\partial \Delta)$.

{\it Step 3.} Now we can prove the Hopf lemma for almost complex
manifolds. Fix a normal field $v$ on $\partial G$ and consider the family
of $J$-holomorphic discs $d_v$ satisfying $d'_0(\partial_x) =
v(d(0))$. The image of such a disc is a real surfaces intesecting
$\partial G$ transversally, so its pullback gives a $\mathcal C^2$-curve in
$\Delta$. Denote by $G_v$ the component of  $\Delta$ defined by the
condition $d_v(G_v) \subset G$. Then every $G_v$ is a domain with
$\mathcal C^2$-boundary in $\C$ and the curvatures of boundaries depend
continuously on $v$. We conclude by applying Step 2 to the composition
$u \circ d_v$ on $G_v$. 

As an application, we obtain the boundary distance preserving property for
biholomorphisms between strictly pseudoconvex domains. 

\begin{proposition}
\label{equiv}
Let $D $ and $D'$ be two  smoothly bounded strictly pseudoconvex
domains in four dimensional almost complex manifolds $(M,J)$ and
$(M',J')$ respectively and let $f:D \rightarrow D'$ be a
$(J,J')$-biholomorphism. Then 
there exists a constant $C > 0$ such that

$$
(1/C) dist(f(z),\partial D') \leq dist(z,\partial D) \leq C dist
(f(z),\partial D').
$$
\end{proposition}
\noindent{\it Proof of Proposition~\ref{equiv}}.
According to the previous section, we may assume that  $D = \{ p:
\rho(p) < 0 \}$ where $\rho$ is a  $J$-plurisubharmonic
function on $D$, strictly $J$-plurisubharmonic in a neighborhood of the
boundary; similarly $D'$ can be defined by means of a function
$\rho'$. Now it suffices
to apply the Hopf lemma to the functions $\rho' \circ f$ and $\rho
\circ f^{-1}$. \qed

\subsection{Boundary continuity of diffeomorphisms}
Using estimates of the
Kobayashi-Royden metric together with the boundary distance preserving
property, we obtain, by means of classical arguments
(see, for instance, K.Diederich-J.E.Fornaess \cite{DiFo}), the following

\begin{proposition}\label{Reg}
Let $D$ and $D'$ be two  smoothly relatively compact strictly pseudoconvex
domains in almost complex manifolds $(M,J)$ and $(M',J')$ respectively. Let
$f:  D  \rightarrow  D'$ be  a smooth
diffeomorphism biholomorphic with respect to $J$ and $J'$. 
Then $f$ extends as a $1/2$-H\"older homeomorphism
between the  closures of $D$ and $D'$.
\end{proposition}

As usual, we denote by $K_{(D,J)}(p,v)$ the value of the
Kobayashi-Royden infinitesimal metric (with respect to the structure $J$) at
a point $p$ and a tangent vector $v$. We begin with the folllowing estimates
of the Kobayashi-Royden infinitesimal metric~:
\begin{lemma}
\label{lowest1}
Let $D$ be a relatively compact strictly pseudoconvex domain in an
almost complex manifold $(M,J)$. Then there
exists a constant $C > 0$ such that

\begin{eqnarray*}
(1/C)\| v \| / dist(p,\partial D)^{1/2} \leq K_{(D,J)}(p,v) 
\leq C\| v \|/dist(p,\partial D)
\end{eqnarray*}for every $p \in D$ and
$v \in T_pM$.
\end{lemma}

\noindent{\it Proof of Lemma~\ref{lowest1}}.
The lower estimate is proved in \cite{GaSu}.
For the upper estimate, it is sufficient to prove the statement near the
boundary. Let $q \in \partial D$. One may suppose that $q = 0$ and  
$J = J_{st} + 0(\vert
z \vert)$. For 
$p \in U \cap D$ consider the ball $p + d(p)\B$, where $d(p)$
is the distance from $p$ to $\partial D$. It follows by
A.Nijenhuis-W.Woolf \cite{NiWo} that there exists  constant $C_1, C_2 > 0$
and a function $d'(p)$
satisfying $C_1d'(p) \leq d(p) \leq C_2d'(p)$ on $D \cap U$, such that for any
complex vector $v$ there exists a $J$ holomorphic map $f:d'(p)\Delta
\rightarrow p + d(p)\B$ such that $f(0) = p$ and
$df_0(e) = v/\| v \|$ ($e$ is the unit vector $1$ in $\C$).
By the decreasing property of the Kobayashi-Royden metric we obtain that 
$K_{(D,J)}(z,v/\| v \| ) \leq
K_{d'(p)\Delta}(0,e)$ which implies the desired estimate.

\vskip 0,1cm
\noindent{\it Proof of Proposition~\ref{Reg}}. For any $p \in D$ and any
tangent vector $v$ at $p$ we have by Lemma~\ref{lowest1}~:

\begin{eqnarray*}
C_1\frac{\| df_p(v) \|}{dist(f(p),\partial D')^{1/2}} \leq 
K_{(D',J')}(f(p),df_p(v)) = K_{(D,J)}(p,v) \leq C_2\frac{\| v
\|}{dist(p,\partial D)}
\end{eqnarray*}
which implies, by Proposition~\ref{equiv}, the estimate
$$\vert\vert\vert df_p \vert\vert\vert \leq C\frac{\| v
\|}{dist(p,\partial D)^{1/2}}.$$
This gives the desired statement. \qed

\vskip 0,2cm
\noindent Proposition~\ref{Reg} allows to reduce the proof of Fefferman's
theorem to a {\it local situation}. 
Indeed, let $p$ be a boundary point of $D$ and
$f(p) = p' \in \partial D'$. It suffices to prove that $f$ extends
smoothly to a neighborhood of $p$ on $\partial D$. Consider
coordinates $z$ and $z'$ defined in small neighborhoods $U$ of $p$ and $U'$ of
$p'$ respectively, with $U' \cap D' = f(D \cap U)$ (this is possible since
$f$ extends as a homeomorphism at $p$). We obtain the following situation.
If $\Gamma = z(\partial D \cap U)$ and $\Gamma' = z'(\partial D' \cap U')$
then the map $z' \circ f \circ z^{-1}$ is defined on $z(D \cap U)$ in $\C^2$,
continuous up to the hypersurface $\Gamma$ with $f(\Gamma) \subset \Gamma'$.
Furthermore the map $z' \circ f \circ z^{-1}$ is a diffeomorphism between
$z(D \cap U)$ and $z'(D' \cap U')$ and the hypersurfaces $\Gamma$ and
$\Gamma'$ are strictly pseudoconvex for the structures $z_*(J)$ and
$(z')_*(J')$ respectively. Finally, we may choose $z$ and $z'$ such that
$z_*(J)$ and $z'_*(J')$ are represented by diagonal matrix functions in the
coordinates $z$ and $z'$. 
As we proved in Lemma~\ref{PP}, $\Gamma$ (resp. $\Gamma'$) is also strictly
$J_{st}$-psdeudoconvex at the origin. We call such
coordinates $z$ (resp. $z'$) {\it canonical coordinates} at $p$
(resp. at $p'$). Using the non-isotropic
dilation as in Section 2.5, we may assume that the norms
$\| z_*(J) - J_{st}\|_{\CC^2}$ and
$\| z'_*(J') - J_{st}\|_{\CC^2}$ are as small as needed.
This localization is crucially used in the sequel and we write $J$ (resp.
$J'$) instead of $z_*(J)$ (resp. $z'_*(J')$); we 
identify $f$ with $z' \circ f \circ z^{-1}$.  

\subsection{Localization and boundary behavior of the tangent map} 

In what follows we will need a more precise information about the boundary
behavior of the tangent map of $f$. Recall that according to the
previous subsection, $D$ and $D'$ are
supposed to be domains in $\C^2$, $\Gamma$ and $\Gamma'$ are open
smooth pieces of their boundaries containing the origin, the almost complex
structure $J$ (resp. $J'$) is defined in a neighborhood of $D$
(resp. $D'$), $f$ is a $(J,J')$ biholomorphism from $D$ to $D'$,
continuous up to $\Gamma$, $f(\Gamma) = \Gamma'$, $f(0) = 0$.
The matrix $J$ (resp. $J'$) is diagonal on $D$ (resp. $D'$).

\vskip 0,1cm
Consider a basis $(\omega_1,\omega_2)$ of $(1,0)$ differential forms
(for the structure $J$) in a neighborhood of the origin. Since $J$ is
diagonal, we may choose $\omega_j = dz^j - B_{j}(z)d\bar z^j$, $j=1,2$.
Denote by $Y=(Y_1,Y_2)$ the corresponding dual basis
of $(1,0)$ vector fields. Then $Y_j = \partial /\partial z^j -
\beta_j(z)\partial/\bar\partial z^j$, $j=1,2$. Here $\beta_j(0) =
\beta_k(0) = 0$. The basis $Y(0)$ simply coincides with the canonical (1,0)
basis of $\C^2$.
In particular $Y_1(0)$ is a basis vector of the holomorphic tangent space
$H^J_0(\partial D)$ and $Y_2(0)$ is ``normal'' to $\partial D$.
Consider now for $t \geq 0$ the translation $\partial D -
t$ of the boundary of $D$ near the origin. Consider, in a neighborhood of the
origin, a $(1,0)$ vector field $X_1$ (for $J$) such that $X_1(0) = Y_1(0)$
and $X_1(z)$ generates the complex tangent space $H^J_z(\partial D - t)$ at
every point $z \in \partial D - t$, $0 \leq t <<1$.
Setting $X_2 = Y_2$, we obtain a basis of vector fields
$X = (X_1,X_2)$ on $D$ (restricting $D$ if necessary). 
Any complex tangent vector $v \in T_z^{(1,0)}(D,J)$ at
point  $z \in D$ admits the unique
decomposition $v = v_t + v_n$ where $v_t = \alpha_1
X_1(z)$ (the tangent component) and $v_n = \alpha_2 X_2(z)$ (the normal
component). Identifying $T_z^{(1,0)}(D,J)$ with $T_zD$ we may
consider the decomposition $v=v_t + v_n$ for $v \in T_z(D)$.
Finally we consider this decomposition for points $z$ in a neighborhood of
the boundary.

 We fix a (1,0) basis vector fields $X$
(resp. $X'$) on $D$ (resp. $D')$ as above.

\begin{proposition}
\label{matrix}
The matrix $A = (A_{kj})_{k,j= 1,2}$ of the differential $df_z$ with respect
to the bases
$X(z)$ and $X'(f(z))$ satisfies the following estimates~:
$A_{11} = O(1)$, $A_{12}= O(dist(z,\partial D)^{-1/2})$,
$A_{21}= O(dist(z,\partial D)^{1/2})$ and $A_{22}= O(1)$.   
\end{proposition}
We begin with the following estimates of the Kobayashi-Royden infinitesimal
pseudometric~:
\begin{lemma}
\label{lowest}
 There exists a positive constant $C$  such that for any $p \in D$
and $v \in T_pD$~:
$$
\frac{1}{C} \left(\frac{\vert v_t \vert}{dist(p,\partial D)^{1/2}} +
\frac{\vert v_n \vert}{dist(p,\partial D)}\right) \leq K_{(D,J)}(p,v)
\leq C \left( \frac{\vert v_t \vert}{dist(p,\partial D)^{1/2}} +
\frac{\vert v_n \vert}{dist(p,\partial D)}\right).
$$
\end{lemma}

\noindent{\it Proof of Lemma~\ref{lowest}}.
The lower estimate is proved in \cite{GaSu}.
For $p \in D$ denote by $d(p)$ the distance from $p$ to $\partial D$.
As in the proof of Lemma~\ref{lowest1}, it follows by
Nijenhuis-Woolf~\cite{NiWo} that there is $r>0$, independent of $p$, such
that for any tangent vector $v$ at $p$ 
satisfying $\omega_2(v) = 0$, there exists a $J$-holomorphic map $f:
r(d(p))^{1/2}\Delta \rightarrow D$ satisfying $df_0(e) = v$. This
implies the upper estimate. \qed

\vskip 0,1cm
\noindent{\it Proof of Proposition~\ref{matrix}}.
Consider the case where $v = v_t$. It follows from
Lemma~\ref{lowest} that~:
$$
\begin{array}{llcll}
\displaystyle \frac{1}{C}\left (
\frac{\| (df_z(v_t))_t \|}{dist(f(z),\partial
D')^{1/2}} + \frac{\|(df_z(v_t))_n
\|}{dist(f(z),\partial D')} \right ) & \leq & 
K_{(D',J')}(f(z),df_z(v_t)) && \\
& = & K_{(D,J)}(z,v_t) & \leq & \displaystyle C
\frac{\| v_t \|}{dist(z,\partial D)^{1/2}}.
\end{array}
$$
This implies that 
$\|(df_z(v_t))_t \| \leq C^{5/2} \| v_t \|$
and
$\vert\vert (df_z(v_t))_n \vert\vert \leq C^{3}dist(z,\partial D)^{1/2}
\| v_t \|$, by the boundary distance
preserving property given in Proposition~\ref{equiv}.
We obtain the estimates for the normal component in a similar way.
\qed

\section{Boundary regularity of a pseudoholomorphic disc attached to a
  totally real manifold}

This section is devoted to one of the main technical steps of our
construction. We prove that a pseudoholomorphic disc attached (in the
sense of the cluster set) to a smooth totally real submanifold in an almost
complex manifold, extends  smoothly up to the boundary. In the case of
the integrable structure, various versions of this statement have been
obtained by several authors. In the almost complex case, similar
assertions have been established by H.Hofer \cite{Ho}, J.-C.Sikorav
\cite{Si}, S.Ivashkovich-V.Shevchishin \cite{IvSh}, E.Chirka \cite{Ch1}
under stronger assumptions on the initial boundary regularity of the disc
(at least the continuity is required). Our proof consists of two
steps. First, we show that a disc extends as a $1/2$-H\"older
continuous map up to the boundary. The proof is based on special
estimates of the Kobayashi-Royden metric in ``Grauert tube'' type
domains. The second step is the reflection principle adapted to the
almost complex category; here we follow the construction of E.Chirka
\cite{Ch1}.

\subsection{H\"older extension of holomorphic discs}
We study the boundary continuity of pseudoholomorphic discs
attached to smooth totally real submanifolds in almost complex manifolds.

Recall that in the case of the integrable structure every smooth totally
real submanifold  $E$ (of maximal dimension) is the zero set of a positive
strictly plurisubharmonic function of class $\CC^2$. This remains true
in the almost complex case. Indeed, we can choose coordinates
$z$ in a neighborhood $U$ of $p \in E$ such that $z(p) = 0$, $z_*(J) =
J_{st} + O(\vert z \vert)$ on $U$ and 
$z(E \cap U) = \{w=(x,y) \in z(U) : r_j(w) = x_j +o(\vert
w \vert) = 0 \}$. The function $\rho = \sum_{j=1}^n r_j^2$ is strictly
$J_{st}$-plurisubharmonic on $z(U)$ and so remains strictly
$z_*(J)$-plurisubharmonic, restricting $U$ if necessary.
Covering $E$ by such neighborhoods, we conclude by mean of the partition of
unity.

Let $f : \Delta \rightarrow (M,J)$ be a $J$-holomorphic disc and let
$\gamma$ be an open arc on the unit circle $\partial \Delta$.
As usual we denote by $C(f,\gamma)$ the cluster set of $f$ on $\gamma$;
this consists of points $p \in M$ such that $p=\lim_{k \rightarrow
\infty}f(\zeta_k)$ for a sequence $(\zeta_k)_k$ in $\Delta$ converging
to a point in $\gamma$.

\begin{proposition}
\label{regth1}
Let $G$ be a relatively compact domain in an almost complex manifold $(M,J)$
and let $\rho$ be a strictly $J$-plurisubharmonic
function of class $\CC^2$ on $\bar{G}$.
Let $f:\Delta \rightarrow G$ be a
$J$-holomorphic disc such that $\rho \circ f \geq 0$ on $\Delta$.
Suppose that $\gamma$ is an open non-empty arc on
$\partial \Delta$ such that the cluster set
$C(f,\gamma)$ is contained in the zero set of $\rho$.
Then $f$ extends as a H\"older 1/2-continuous map on $\Delta \cup
\gamma$.
\end{proposition}

We begin the proof by the following well-known assertion
(see, for instance, \cite{BeL}).

\begin{lemma}
\label{dlem3.1}
Let $\phi$ be a positive subharmonic function in $\Delta$ such that
the measures $\mu_r(e^{i\theta}) := \phi(re^{i\theta})d\theta$ converge in
the weak-star topology to
a measure $\mu$ on $\partial \Delta$ as $r \rightarrow 1$. Suppose that
$\mu$ vanishes on an open arc $\gamma \subset \partial \Delta$. Then for
every compact subset $K \subset \Delta \cup \gamma$ there exists a constant
$C>0$ such that
$\phi(\zeta) \leq C(1 - \vert \zeta \vert)$ for any
$\zeta \in K \cup \Delta$.
\end{lemma}
 
Now fix a point $a \in \gamma$, a constant $\delta > 0$ small
enough so that the intersection $\gamma \cap (a + \delta
\bar\Delta )$ is compact in $\gamma$; we denote by
$\Omega_{\delta}$ the intersection $\Delta \cap (a + \delta\Delta
)$. By Lemma~\ref{dlem3.1}, there exists a constant $C > 0$ such that,
for any $\zeta$ in $\Omega_{\delta}$, we have

\begin{eqnarray}
\label{dd4}
\rho \circ f(\zeta) \leq C (1 - \vert \zeta \vert ).
\end{eqnarray}

Let $(\zeta_k)_k$ be a sequence of points in $\Delta$ converging to $a$
with $\lim_{k \rightarrow \infty}f(\zeta_k) = p$.
By assumption, the function $\rho$ is strictly $J$-plurisubharmonic in a
neighborhood $U$ of $p$; hence there is a constant
$\varepsilon > 0$ such that the function $\rho - \varepsilon \vert z
\vert^2$ is $J$-plurisubharmonic on $U$.

\begin{lemma}
\label{dlem3.2}
There exists a constant $A > 0$ with the following property~: If
$\zeta$ is an arbitrary point of $\Omega_{\delta/2}$ such that
$f(\zeta)$ is in $G \cap z^{-1}(\B)$, then
$\vert \vert \vert df_\zeta \vert \vert \vert
\leq A(1- \vert \zeta \vert)^{-1/2}$.
\end{lemma}
\noindent{\it Proof of Lemma~\ref{dlem3.2}.}
Set $d = 1 - \vert \zeta \vert$; then the disc $\zeta +
d\Delta$ is contained in $\Omega_{\delta}$. Define the domain $G_d =
\{ w \in G: \rho(w) < 2Cd \}$. Then it follows by (\ref{dd4}) that the
image $f(\zeta + d\Delta)$ is contained in $G_d$, where the
$J$-plurisubharmonic function $u_d = \rho - 2Cd$ is negative. 
Moreover we have the following lower estimates on the Kobayashi-Royden
infinitesimal pseudometric (a rather technical proof is given in the
Appendix)~:

\begin{proposition}
\label{lowest0,5}
Let $D$ be a domain in an almost complex manifold $(M,J)$, let $p \in
\bar D$, let $U$ be a neighborhood of $p$ in $M$ (not necessarily
contained in $D$) and let $z: U \rightarrow  B$ be a
normalized coordinate diffeomorphism introduced in Lemma~\ref{suplem1}.
Let $u$ be a $\CC^2$ function on
$D$, negative and $J$-plurisubharmonic on $D$. We assume that $-L \leq
u < 0$ on $D \cap U$ and that $u - c \vert z \vert^2$ is
$J$-plurisubharmonic
on $D \cap U$, where $c$ and $L$ are positive constants. Then there
exists a neighborhood $U' \subset U$ of $p$ depending on $c$ and
$\| J \|_{C^2(U)}$, a positive constant $c'$, depending only on $c$
and $L$, such that we have the following estimate:

$$K_{(D,J)}(q,v) \geq c'\| v \| /\vert u(q)
\vert^{1/2}$$
for every $q \in D \cap U'$ and every $v \in T_qM$.
\end{proposition}

Hence there exists a positive constant $M$
(independent of $d$) such that $K_{(G_d,J)}(w,\eta) \geq M \vert \eta
\vert \vert u_d(w) \vert^{-1/2}$, for any $w$ in $G \cap z^{-1}(\B)$ and
any $\eta \in T_{w}\Omega$. On another hand, we have $K_{\zeta +
d\Delta}(\zeta,\tau ) = \vert \tau \vert /d$ for any $\tau$ in
$T_{\zeta}\Delta$ indentified with $\C$. By the decreasing property
of the Kobayashi-Royden metric, for any $\tau$ we have 

\begin{eqnarray*}
M \| df_\zeta(\tau)\| \ \vert u_d(f(\zeta))\vert^{-1/2} \leq
K_{(G_d,J)}(f(\zeta),df_\zeta(\tau)) \leq K_{\zeta + d\Delta}(\zeta,\tau) =
\vert \tau \vert/d.
\end{eqnarray*}
Therefore, $\vert \vert \vert df_\zeta\vert\vert\vert \leq M^{-1}\vert
u_d(f(\zeta))\vert^{1/2}/d$. As $-2Cd \leq u_d(f(\zeta)) < 0$, this
implies the desired statement in Lemma~\ref{dlem3.2}
with $A = M^{-1}(2C)^{1/2}$. \qed

\vskip 0,1cm
\noindent{\it Proof of Proposition~\ref{regth1}}.
Lemma~\ref{dlem3.2} implies that $f$ extends  as a 1/2-H\"older map
to a neighborhood of the point $a$
in view of an integration argument inspired by the classical
Hardy-Littlewood theorem. 
This proves Proposition~\ref{regth1}. \qed

\subsection{Reflection principle and regularity of analytic discs}

In the previous subsection we proved that a $J$-holomorphic disc attached to
a smooth totally real submanifold is $1/2$-H\"olderian up to the
boundary. This allows to use the reflection principle for pseudoholomorphic
curves. Similar ideas have been used by E.Chirka \cite{Ch1}
and S.Ivashkovich-V.Shevchischin \cite{IvSh}.
For reader's convenience we present the argument due to E.~Chirka.
\begin{proposition}
Let $E$ be an $n$-dimensional smooth totally real submanifold  in an
almost complex manifold $(M,J)$.  For any $p \in E$ there exists a
neighborhood $U$ of $p$ and a smooth coordinate diffeomorphism 
$z:U \rightarrow \B$ such that $z(E) = \R^n$ and $z_*(J) _{\vert
\R^n} = J_{st}$. Moreover, the condition of $z_*(J)$-holomorphy for a
disc $f$ may be written in the form
$$\bar\partial f + A(f)\overline{\partial f} = 0$$
where the smooth matrix function $A(z)$ vanishes with infinite order
on $\mathbb R^n$.
\end{proposition}  

\proof After a complex linear change of coordinates we may assume that
$J = J_{st} + O(\vert z \vert)$ and $E$ is given by $x + ih(x)$ where
$x \in \R^n$ and $dh(0) = 0$. If $\Phi$ is the local diffeomorphism
$x \mapsto x$, $y \mapsto y - h(x)$ then $\Phi(E) = \R^n$ and the direct
image of $J$ by $\Phi$, still denoted by $J$, keeps the form $J_{st} +
O(\vert z \vert)$. Then $J$ has a basis of $(1,0)$-forms given in the
coordinates $z$ by $dz^j + \sum_k a_{jk}d\bar z^k$; using the
matrix notation we write it in the form $\omega = dz + A(z)d\bar z$ where
the matrix function $A(z)$ vanishes at the origin. Writing
$\omega = (I + A)dx + i(I - A)dy$ where $I$ denotes the identity
matrix, we can take as a basis of $(1,0)$ forms~: $\omega' = dx +
i(I +A)^{-1}(I - A)dy = dx + iBdy$. Here the matrix function $B$ satisfies
$B(0) = I$. Since $B$ is smooth, its restriction $B_{\vert \R^n}$ on $\R^n$
admits a smooth extension $\hat B$ on the unit ball such that
$\hat B - B_{\vert \R^n} = O(\vert y \vert^k)$ for any positive integer $k$.
Consider the diffeomorphism $z^* = x + i\hat B(z) y$.
In the $z^*$-coordinates the submanifold $E$ still coincides with $\R^n$
and $\omega' = dx + iBdy = dz^* + i(B - \hat B)dy - i(d\hat B)y = dz^* +
\alpha$, where the coefficients of the form $\alpha$ vanish with infinite
order on $\R^n$. Therefore there is a basis of $(1,0)$-forms
(with respect to the image of $J$ under the coordinate diffeomorphism
$z \mapsto z^*$) of the form $dz^* + A(z^*)d\bar z^*$,
where $A$ vanishes with infinite order on $\R^n$. \qed

\vskip 0,1cm
Now we are able to prove the main result of this section. We denote by
$\Delta^+$ the upper half disc $\Delta^+ = \{\zeta \in \C : Im(\zeta) > 0\}$.

\begin{proposition}\label{prop-reg}
Let $E$ be a smooth totally real $n$-dimensional submanifold in a real
$2n$-dimensional almost complex manifold $(M,J)$ and let
$f:\Delta^+ \rightarrow M$ be a
$J$-holomorphic map. If the cluster set $C(f,]-1,1[)$ is
(compactly) contained in $E$, then
$f$ is of class $\CC^{\infty}$ on $\Delta^+ \cup ]-1,1[$. 
\end{proposition}

\proof We know that $f$ is $1/2$-H\"older continuous on $\Delta^+ \cup
]-1,1[$ by Proposition~\ref{regth1} (see the beginning of Subsection~4.1).
Moreover, it is easy to see that the H\"older constant
depends on $\|f\|_\infty$, on the Levi
form of the positive strictly $J$-plurisubharmonic defining function of
$E$ and on the $\CC^2$-norm of $J$.
Fix a point $a \in ]-1,1[$. We may assume that $g(a)  := z \circ f(a) =
0$. It is enough to show that the map $g$ is of
class $\CC^{\infty}$ on $\Delta^+ \cap ]-1,1[$ in a neighborhood of
$a$. Consider the map $\hat g$ equal to $g$ on $\Delta^+ \cup
]-1,1[$
and defined by $\hat g(\zeta) = \overline{g(\bar \zeta)}$ for
$\zeta \in \Delta^-$. Then $\hat g$ is continuous on $\Delta$. We
write
the $z_*(J)$-holomorphy condition for $g$ on $\Delta^+$ in the form
$\bar \partial g + q(g)\overline{\partial g} = 0$
where $q$ is a smooth matrix function satisfying
$\| q \|_{\mathcal C^k} < < 1$ for any $k$. Then for $\zeta
\in \Delta^-$ we have $\bar\partial\hat g(\zeta) + 
\overline{q(g(\bar\zeta))}\ \overline{\partial\hat g(\zeta)} = 0$.
This means that $\hat g$ satisfies on $\Delta$ an elliptic equation of the
form $\bar \partial \hat g + \phi(\cdot)\overline{\partial \hat g} = 0$
where $\phi$ is defined by $\phi(\zeta) = q(g(\zeta))$ for $\zeta \in
\Delta^+ \cup ]-1,1[$ and $\phi(\zeta) =
\overline{q(g(\bar\zeta))}$ for $\zeta \in \Delta^-$.
Since $q$ vanishes on $\R^n$ with
infinite order and $g$ is $1/2$-H\"older continuous up to $]-1,1[$, the
matrix function $\phi$ is smooth on $\Delta$.
The result now follows by the well-known results on 
elliptic regularity (see, for instance, \cite{NiWo,Si}) since $\hat g$
is necessarily $\CC^{\infty}$-smooth on $\Delta$. \qed

\vskip 0,1cm
Moreover, it follows from these results that for any $k$ and for any $r > 0$
there exists a constant such that
$\| \hat g \|_{\mathcal C^k((1-r)\Delta)} \leq
C \| \hat g \|_{\mathcal C^0((1-r)\Delta)}$.
This estimate will be used in the following section to study
the boundary regularity of pseudoholomorphic maps. Since the initial
disc $f$ is obtained from $\hat g$ by a diffeomorphism depending on
$E$ and $J$ only, we obtain the
following quantitative version of the previous statement.

\begin{proposition}
Let $E$ be a totally real $n$-dimensional submanifold in an
(complex) $n$-dimensional almost complex manifold $(M,J)$
and let $f:\Delta^+ \rightarrow \C^n$ be a
$J$-holomorphic map. Assume that the cluster set $C(f,]-1,1[)$ is
compactly contained in $E$. Then given $r > 2$ there exists a constant
$C >0$ depending on $\| f \|_{\infty}$ and on the $\CC^r$ norm of the
defining functions of $E$ such that 
\begin{eqnarray}
\label{uniformity}
\| f \|_{\mathcal C^r(\Delta^+ \cup ]-1,1[)} \leq C 
\| J \|_{\mathcal C^r}.
\end{eqnarray}
\end{proposition}

In the next Section we apply these results to the study of the bounadary
regularity of pseudoholomorphic maps of wedges with totally real edges.
                              
\section{Behavior of pseudoholomorphic maps near totally real submanifolds}

Let $\Omega$ be a domain in an almost complex manifold $(M,J)$ and $E 
\subset \Omega$ be a smooth $n$-dimensional
totally real submanifold defined as the set of common zeros of the
functions $r_j$, $j=1,...,n$ smooth on $\Omega$. We suppose that
$\bar\partial_J r_1 \wedge ...\wedge \bar\partial_J r_n \neq
0$ on $\Omega$. Consider the ``wedge''
$W(\Omega,E)=\{ z \in \Omega: r_j(z) < 0, j= 1,...,n \}$ with ``edge''
$E$. For $\delta > 0$ we denote by $W_{\delta}(\Omega,E)$ the
``shrinked'' wedge 
$\{ z \in \Omega : r_j(z) - \delta \sum_{k \neq j} r_k <
0, j = 1,..., n \}$.
The main goal of this Section is to prove the following 

\begin{proposition}
\label{wedges}
Let $W(\Omega,E)$ be a wedge in $\Omega \subset (M,J)$ with a totally real
n-dimensional edge $E$ of class $\CC^{\infty}$ and let $f:W(\Omega,E)
\rightarrow (M',J')$ be a  $(J,J')$-holomorphic map. Suppose that the
cluster set $C(f,E)$ is (compactly) contained in a
$\mathcal C^\infty$ totally real submanifold $E'$ of $M'$.
Then for any $\delta > 0$ the map $f$ extends to
$W_{\delta}(\Omega,E) \cup E$ as a $\CC^{\infty}$-map. 
\end{proposition}

In Section~4 we established this statement for a single
$J$-holomorphic disc. The general case also 
relies on the ellipticity of the
$\bar\partial$-operator.  It requires an additional
 technique of attaching pseudoholomorphic discs to a totally real
manifold which could be of independent interest.

\subsection{Almost complex perturbation of discs}

In this subsection we attach Bishop's discs to a totally real
submanifold in an almost complex manifold.
The following statement is an almost complex
analogue of the well-known Pinchuk's construction \cite{Pi74} of a
family of holomorphic discs attached to a totally real manifold.

\begin{lemma}\label{lem-discs}
For any $\delta > 0$ there exists a family of
$J$-holomorphic discs $h(\tau,t) = h_t(\tau)$ smoothly depending on the
parameter $t \in  \R^{2n}$ such that $h_t(\partial \Delta^+) \subset E$,
$h_t(\Delta) \subset W(\Omega,E)$, $W_{\delta}(\Omega,E) \subset \cup_t
h_t(\Delta)$ and $C_1(1 - \vert \tau \vert) \leq dist (h_t(\tau),E)
\leq C_2 ( 1- \vert \tau \vert)$ for any $t$ and any $\tau \in \Delta^+$,
with constants $C_j > 0$ {\it independent} of $t$.
\end{lemma}
For $\alpha > 1$, noninteger, we denote by $\mathcal C^\alpha(\bar \Delta)$
the Banach space of functions of class $\mathcal C^\alpha$ on $\bar{\Delta}$
and by $\mathcal A^\alpha$ the Banach subspace of
$\mathcal C^\alpha(\bar \Delta)$ of functions holomorphic on $\Delta$.

First we consider the situation where $E=\{r:=(r_1,\dots,r_n)=0\}$
is a smooth totally real submanifold in $\C^n$.
Let $J_{\lambda}$ be an almost complex deformation of the standard
structure $J_{st}$ that is a one-parameter family of almost complex
structures so that $J_0 = J_{st}$. 
We recall that for $\lambda$ small enough the
$(J_{st},J_{\lambda} )$-holomorphy condition for a map $f:\Delta
\rightarrow \C^n$ may be written in the form

\begin{eqnarray}\label{equa0}
\bar\partial_{J_{\lambda}} f = \bar\partial f +
q(\lambda,f)\overline{\partial f} = 0
\end{eqnarray}
where $q$ is a smooth matrix satisfying $q(0,\cdot) \equiv 0$, uniquely
determined by $J_\lambda$ (\cite{Si}).

A disc $f \in  (\mathcal C^\alpha(\bar{\Delta}))^n$ is attached
to $E$ and is $J_\lambda$-holomorphic if and only if it
satisfies the following nonlinear boundary Riemann-Hilbert type  problem~:
$$
\left\{
\begin{array}{lll}
r(f(\zeta)) = 0,& & \zeta \in \partial \Delta\\
\bar{\partial}_{J_\lambda}f(\zeta) = 0,& & \zeta \in \Delta.
\end{array}
\right.
$$
Let $f^0 \in \mathcal (\mathcal A^\alpha)^n$ be a disc attached to
$E$ and let $\mathcal U$ be a neighborhood of $(f^0,0)$ in the 
space $(\mathcal C^\alpha(\bar{\Delta}))^n \times \R$.
Given $(f,\lambda)$ in $\mathcal U$ define the maps 
$v_{f}: \zeta \in \partial 
\Delta \mapsto r(f(\zeta))$
and

\begin{eqnarray*}
& &u : \mathcal{U}  \rightarrow  (\mathcal C^\alpha(\partial \Delta))^n
\times \mathcal C^{\alpha-1}(\Delta)\\
  & & (f,\lambda)  \mapsto  (v_{f}, 
\bar{\partial}_{J_\lambda}f).
\end{eqnarray*}

Denote by $X$ the Banach space $(\mathcal C^\alpha(\bar \Delta))^n$.
Since $r$ is of class $\CC^{\infty}$, 
the map
$u$ is smooth and the tangent map $D_Xu(f^0,0)$ (we consider 
the derivative
with respect to the space $X$) is a linear map from $X$ to 
$(\CC^\alpha(\partial \Delta))^n \times \CC^{\alpha-1}(\Delta)$,
defined for every $h \in X$ by 
$$
D_Xu(f^0,0)(h) = 
\left(
\begin{array}{ll}
0 & 2 Re [G h] \\
0 & \bar\partial_{J_0} h
\end{array}
\right),$$
where for $\zeta \in \partial \Delta$
$$
G(\zeta) = \left(
\begin{array}{lll}
\frac{\partial r_1}{\partial z^1}(f^0(\zeta))  &\cdots&\frac{\partial
r_1}{\partial z^n}(f^0(\zeta))\\
\cdots&\cdots&\cdots\\
\frac{\partial r_n}{\partial z^1}(f^0(\zeta))& \cdots&\frac{\partial r_n}
{\partial z^n}(f^0(\zeta))
\end{array}
\right)$$
(see \cite{gl94}).
\begin{lemma}\label{tthh}
Assume that for some $\alpha > 1$
the linear map from $(\mathcal A^{\alpha})^n$ to
$(\mathcal C^{\alpha-1}(\Delta))^n$
given by $h \mapsto 2 Re [G h]$
is surjective and has a $d$-dimensional kernel.
Then there exist $\delta_0, \lambda_0 >0$ such that for every
$0 \leq \lambda \leq \lambda_0$,
the set of $J_\lambda$-holomorphic discs $f$ attached to $E$
and such that $\| f -f^0 \|_{\alpha} \leq \delta_0$ forms
a smooth $d$-dimensional 
submanifold 
$\mathcal A_{\lambda}$ in the Banach space 
$(C^\alpha(\bar{\Delta}))^n$.
\end{lemma} 

\noindent{\it Proof of Lemma~\ref{tthh}.}
According to the implicit function Theorem, 
the proof of Lemma~\ref{tthh} reduces to the proof of the 
surjectivity of $D_Xu$. 
It follows by classical
one-variable results on the resolution of the
$\bar\partial$-problem in the unit disc that the linear map from
$X$ to $\CC^{\alpha-1}(\Delta)$ given by 
$h \mapsto \bar \partial h$
is surjective. More precisely, given $g \in \CC^{\alpha -1}(\Delta)$
consider the Cauchy transform

$$T_{\Delta}(g) : \tau \in \partial \Delta \mapsto   
\frac{1}{2i\pi}
\int\int_{\Delta} \frac{g(\zeta)}{\zeta - \tau}d\zeta \wedge d\bar{\zeta}.$$

For every function $g \in 
C^{\alpha-1}(\Delta)$ the solutions $h \in X$ of the equation
$\bar\partial h = g$ have the form $h = h_0 + T_{\Delta}(g)$
where $h_0$  is an arbitrary function in $({\mathcal A}^{\alpha})^n$. 
Consider the equation 

\begin{equation}\label{EQU}
D_Xu(f^0,0)(h) = \left(
\begin{array}{lll}
0 & & g_1 \\
0 & & g_2
\end{array}
\right)
\end{equation}
where $(g_1,g_2)$ is a vector-valued function with components 
$g_1 \in \CC^{\alpha-1}(\partial \Delta)$ and
$g_2 \in \CC^{\alpha-1}(\Delta)$.
Solving the $\bar\partial$-equation for the second component, we reduce 
equation~(\ref{EQU}) to    
$$
2 Re [G(\zeta) h_0(\zeta)] = g_1 - 2 Re [G(\zeta) T_{\Delta}(g_2)(\zeta)]
$$
with respect to $h_0 \in (\mathcal A^{\alpha})^n$. 
The surjectivity of the map $ h_0 \mapsto 2 Re [G h_0]$ gives the result.
\qed

\vskip 0,1cm
\noindent{\it Proof of Lemma~\ref{lem-discs}.} We proceed in three steps.
{\it Step 1. Filling the polydisc.} Consider the $n$-dimensional real 
torus $\mathbb T^n = \partial \Delta \times
...\times \partial \Delta$ in $\C^n$ and the
linear disc $f^0(\zeta) = (\zeta,...,\zeta)$, $\zeta \in \Delta$
attached to $\mathbb T^n$. 
In that case, a disc $h^0$ is in the kernel of
$h \mapsto 2 Re [G h]$ if and only if every component $h^0_k$ of $h^0$
satisfies on $\partial \Delta$ the condition $h^0_k +
\zeta^2\overline{h^0_k} = 0$. Considering the Fourier expansion
of $h_k$ on $\partial \Delta$ (recall that $h_k$ is holomorphic on
$\Delta$) and identifying the coefficients, we obtain that the map  
$h \mapsto 2 Re [G h]$ from $({\mathcal A}^{\alpha})^n$ to
$(C^{\alpha - 1}(\Delta))^n$ is surjective and has a $3n$-dimensional
kernel.
By Lemma~\ref{tthh} if $J_\lambda$ is an almost complex
structure close enough to $J_{st}$ in a neighborhood of the closure
of the polydisc $\Delta^n$, there is a $3n$-parameters family of
$J_\lambda$-holomorphic discs attached to $\mathbb T^n$. These
$J_{\lambda}$-holomorphic discs fill the intersection of 
a sufficiently small neighborhood of the point $(1,...,1)$
with $\Delta^n$.

{\it Step 2. Isotropic dilations.} Consider a smooth totally real submanifold
$E$ in an almost complex manifold $(M,J)$. Fixing local coordinates,
we may assume that $E$ is a submanifold in a neighborhood of the
origin in $\C^n$, $J = J_{st} + 0(\vert z \vert)$ and $E$ is defined
by the equations $y = \phi(x)$, where $\nabla \phi(0) = 0$. For
every $\varepsilon > 0$, consider the isotropic dilations
$\Lambda_{\varepsilon}: z \mapsto z' = \varepsilon^{-1}z$. Then
$J_{\varepsilon}:= \Lambda_{\varepsilon}(J) \rightarrow J_{st}$ as
$\varepsilon \rightarrow 0$. In the $z'$-coordinates $E$ is defined
by the equations $y' = \psi(x',\varepsilon):=
\varepsilon^{-1}\phi(\varepsilon x')$ and $\psi \rightarrow 0$
as $\varepsilon \rightarrow 0$. Consider the local diffeomorphism
$\Phi_{\varepsilon}: z' = x' +iy' \mapsto
z''= x' +i(y'-\psi(x',\varepsilon))$. Then in new coordinates (we omit
the primes) $E$ coincides with a neighborhood of the origin in $\R^n = \{
y = 0  \}$ and $\hat J_{\varepsilon}: =
(\Phi_{\varepsilon})_*(J_{\varepsilon}) \rightarrow J_{st}$ as
$\varepsilon \rightarrow 0$. Furthermore, applying a
fractional-linear transformation of $\C^n$, biholomorphic with
respect to $J_{st}$, we may assume that $E$ is a
neighborhood of the point $(1,...,1)$ on the torus $\mathbb T^n$ and the
almost complex structure $J_{\varepsilon}$ is a small deformation of
the standard structure. By Step 1, we may fill a neigborhood of
the point $(1,...,1)$ in the polydisc $\Delta^n$ by
$J_{\varepsilon}$-holomorphic discs (for $\varepsilon$ small enough) which
are small perturbations of the disc $\zeta \mapsto (\zeta,...,\zeta)$.
Returning to the initial coordinates, we obtain a family of $J$-holomorphic
discs attached to $E$ along a fixed arc (say, the upper semi-circle $\partial
\Delta^+$) and filling the intersection of a neighborhhod of the origin with
the wedge $\{y - \phi(x) < 0\}$. 

{\it Step3.} Let now $W(\Omega,E) = \{ r_j < 0 , j=1,...,n\}$ be a
wedge with edge $E$; we assume that $0 \in E$ and $J(0) =
J_{st}$. We may assume that $E=\{y = \phi(x)\}$, $\nabla \phi(0) = 0$,
since the linear part of every $r_j$ at the origin is equal to $y_j$. So
shrinking $\Omega$ if necessary, we obtain that for any $\delta > 0$
the wedge $W_{\delta}(\Omega, E) = \{z \in \Omega: r_j(z) -
\delta\sum_{k \neq j} r_k(z) < 0 , j=1,...,n \}$ is contained in the
wedge $\{z \in \Omega: y - \phi(x) < 0 \}$. By Step 2 there is a family of
$J$-holomorphic discs attached to $E$ along the upper semi-circle and
filling the wedge $W_{\delta}(\Omega,E)$. These discs are  smooth up to the
boundary and smoothly depend on the parameters. \qed  

\subsection{Uniform estimates of derivatives.} 
Now we prove Proposition~\ref{wedges}.
Let $(h_t)_t$ be the family of $J$-holomorphic discs, smoothly depending on
the parameter $t \in \R^{2n}$, defined in Lemma~\ref{lem-discs}.
It follows from Lemma~\ref{dlem3.2}, applied to the holomorphic disc $f \circ
h_t$, uniformly with respect to $t$, that there is a constant $C$
such that $\vert \vert \vert df(z) \vert \vert \vert \leq C dist(z,E)^{-1/2}$
for any $z \in W_{\delta}(\Omega,E)$.
This implies that $f$ extends as a
H\"older $1/2$-continuous map on $W_{\delta}(\Omega,E) \cup E$. 

It follows now from Proposition~\ref{prop-reg} that every
composition $f \circ h_t$ is smooth up to $\partial \Delta^+$. Moreover,
since $f$ is continuous up to $E$, the estimate (\ref{uniformity}) shows
that in our case the $\CC^k$ norm of the discs $f \circ h_t$ are uniformly
bounded, for any $k$. Recall the separate smoothness principle
(Proposition 3.1, \cite{Tu94}):

\begin{proposition}\label{separate}
Let $F_j$, $1 \leq j \leq n$, be $\CC^{\alpha}$ ($\alpha > 1$
 noninteger)
 smooth foliations in a domain $\Omega \subset \R^n$ such
that for every point $p \in \Omega$ the tangent vectors to the curves
$\gamma_j \in F_j$ passing through $p$ are linearly independent. Let
$f$ be a function on $\Omega$ such that the restrictions $f
_{\vert{\gamma_j}}$, $1 \leq j \leq n$, are of class
$\CC^{\alpha-1}$ and are uniformly bounded in the $\CC^{\alpha-1}$
norm. Then $f$ is of class $\CC^{\alpha-1}$.
\end{proposition}
Using Lemma~\ref{lem-discs} we construct $n$ transversal
foliations of $E$ by boundaries of Bishop's discs. Since the restriction
of $f$ on every such curve satisfies the hypothesis of
Proposition~\ref{separate}, $f$ is smooth up to $E$. This proves
Proposition~\ref{wedges}. \qed

\section{Lifts of biholomorphisms to the cotangent bundle}

We first recall the notion of conormal bundle of a real submanifold in $\C^n$
with the standard structure (\cite{tu01}).
Let $T^*(\C^n)$  be the real cotangent bundle of $\C^n$, identified with
the bundle $T^*_{(1,0)}(\C^n)$ of complex (1,0)
forms and let $\pi: T^*(\C^n) \rightarrow
\C^n$ be the natural projection. 
In the canonical complex coordinates $(z,t)$ on $T^*_{(1,0)}(\C^n)$
an element of the fiber at a point $z \in \C^n$ is a (1,0) 
form $\omega =\sum_j t_j dz^j$. Let $N$ be a real smooth generic submanifold
in $\C^n$. The conormal bundle $\Sigma(N)$ of $N$ is
a real subbundle of $T^*_{(1,0)}(\C^n)$ defined by the condition
$\Sigma(N) = \{ \phi \in T^*_{(1,0)}(\C^n): Re
\,\phi \vert T_z^{(1,0)}(N) = 0 , z \in N \}$. By $T_z^{(1,0)}(N)$ we
mean here the real tangent space of $N$ at $z$ considered as a
subspace in $T^{1,0}_z(\C^n)$ after the canonical identification of
the tangent bundles $T(\C^n)$ and $T^{(1,0)}(\C^n)$.

Let $\rho_1,\dots,\rho_d$ be local defining functions of $N$.
Then the forms $\partial \rho_1,\dots, \partial \rho_d$ form
a basis in $\Sigma_{z}(N)$ and every section $\phi$ of the bundle
$\Sigma(N)$ has the form $\phi = \sum_{j=1}^d c_j \partial \rho_j$, 
$c_1,\dots,c_d \in
\R$. We will use the following statement (see~\cite{tu01})~:
\begin{lemma}\label{lemma}
Let $\Gamma$ be a $\CC^2$  real hypersurface in $\C^n$. 
The conormal bundle $\Sigma(\Gamma)$ (except the zero section) is 
a totally real
submanifold of dimension $2n$ in $T^*_{(1,0)}(\C^n)$ if and 
only if the Levi form of $\Gamma$ is nondegenerate.
\end{lemma} 

The conormal bundle notion can be easily extended to the case of an almost
complex manifold.  Let $i: T^*(M) \rightarrow  T^*_{(1,0)}(M,J) $ be the
canonical identification. Let $D$ be a smoothly relatively compact
domain in $M$ with boundary $\Gamma$. The conormal bundle
$\Sigma_J(\Gamma)$ of $\Gamma$ is the real subbundle of
$T^*_{(1,0)}(M,J)$ defined by
$\Sigma_J(\Gamma) = \{ \phi \in T^*_{(1,0)}(M,J): Re \,\phi \vert
T_{z}^J(\Gamma) = 0, z \in \Gamma \}$. As above, by $T_{z}^J(\Gamma)$
we mean the real tangent space  of $\Gamma$ at $z$ viewed as a real
subspace of the (1,0) (with respect to $J$) tangent space of $M$
at $z$.

This notion is invariant with respect to
biholomorphisms. More precisely, if $f:(D,J) \rightarrow (D',J')$
is a biholomorphic map $\CC^1$-smooth up to $\partial D$, then its {\it
cotangent map} defined by $\tilde f:= (f,{}^tdf^{-1})$ is continuous
up to $\Sigma_J(\partial D)$ and $\tilde{f}(\Sigma_J(\partial
D)) = \Sigma_{J'}(\partial D')$. To apply the results of the previous
sections, we define an almost complex structure
on the cotangent bundle $T^*(M)$ of an almost complex manifold such that the
cotangent map of a biholomorphism is biholomorphic with respect to this
structure. For reader's convenience we recall the explicit construction of
this almost complex structure $\tilde{J}$ (ie. the proof of the following
Proposition), following~\cite{YI}, in Appendix~2. 

\begin{proposition}\label{prop-lift}
Let $(M,J)$ be an almost complex manifold. There exists an almost complex
structure $\tilde{J}$ on $T^*M$ with the following properties~:

$(i)$ If $f$ is a biholomorphism between $(M,J)$ and $(M',J')$
then the cotangent map $\tilde f$ is a biholomorphism between
$(T^*M,\tilde J)$ and $(T^*M',\tilde J')$.

$(ii)$ If $(J_\varepsilon)_\varepsilon$ is a small deformation of the
standard structure on $\C^n$ then
$\tilde{J}_\varepsilon \rightarrow J_{st}$ as $\varepsilon \rightarrow 0$
in the $\CC^k$-norm on $T^*\C^n$ (for any $k$).
\end{proposition}

Consider now a smooth relatively compact strictly pseudoconvex domain $D$
in an almost complex manifold $(M,J)$ of real dimension four.
We have the following 
\begin{lemma}\label{totreal}
The conormal bundle of $\partial D$ (outside the zero section)
is a totally real submanifold in $(T^*(M),\tilde J)$.
\end{lemma}

\noindent{\sl Proof of Lemma~\ref{totreal}}. According to Section~2 we may
choose local coordinates near a boundary point $p$ such that $\partial D$
is given by the equation $Re z^2 + Re K(z) + H(z) + o(\vert z \vert^2)
= 0$ and $H(z^1,0)$ is a positive definite hermitian form on $\C$; in
these coordinates the matrix $J$ is diagonal and $J(0) =
J_{st}$. After the non-isotropic dilation $(z^1,z^2) \mapsto
(\varepsilon^{-1/2}z^1,\varepsilon^{-1}z^2)$ the hypersurface is
defined by the equation 
$Re z^2 + \varepsilon^{-1}K(\varepsilon^{1/2}z^1,\varepsilon z^2)
+ H(\varepsilon^{1/2}z^1,\varepsilon z^2) +\varepsilon^{-1}
o(\vert (\varepsilon^{1/2}z^1,\varepsilon z^2)\vert^2) = 0$
and the dilated structure, denoted by $J_\varepsilon$, tends
(together with all derivatives of any order) 
to the standard structure. The hypersurface $\partial D$ tends to the
strictly $J_{st}$-pseudoconvex hypersurface $\Gamma_0 = \{ Re z^2 + K(z^1,0)
+ H(z^1,0) = 0 \}$. It follows from Proposition~\ref{prop-lift} $(ii)$
that $\tilde J_{\varepsilon}$ tends to the standard complex structure
on $T^*(\C^n)$ as $\varepsilon \rightarrow 0$. Since the conormal
bundle of $\Gamma_0$ (with respect to $J_{st}$) is totally real with
respect to the standard structure on $T^*(\C^n)$, the
same holds for $\Sigma_J(\partial D)$ in a small neighborhood of $p$,
by continuity (for $\varepsilon$ small enough). \qed

\vskip 0,1cm
If $f: (D,J) \rightarrow (D',J')$ is a biholomorphism between two
strictly pseudoconvex domains, of class $\CC^1$ on $\bar D$, then
its cotangent lift extends continuously on $\Sigma_J(\partial D)$ and
$f(\Sigma_J(\partial D) \subset \Sigma_{J'}(\partial D')$.
In view of Proposition~\ref{wedges} this proves Fefferman's theorem
under the additional assumption of $\CC^1$-smoothness of $f$ up to the
boundary~:

\begin{proposition}
Let $D$ and $D'$ be smooth relatively compact strictly pseudoconvex domains
in (real) four dimensional almost complex manifolds $(M,J)$ and
$(M',J')$. Consider a $\CC^{\infty}$- diffeomorphism  $f:(D,J)
\rightarrow (D',J')$  which is a  $\CC^1$-diffeomorphism between
$\bar D$ and $\bar D'$. Suppose that the direct image 
$ f_*(J)$ extends $\CC^{\infty}$-smoothly on $\bar D'$. Then $f$
is a $\CC^{\infty}$-diffeomorphism between $\bar D$
and $\bar D'$.
\end{proposition}

This statement just follows by the reflection
principle of Sections~4 and 5.
In order to get rid of the $\CC^1$-assumption we will use the
estimates on the Kobayashi-Royden metric and the scaling method.

\section{Scaling on almost complex manifolds}

Our goal now is to prove Fefferman's mapping theorem without the
assumption of $\CC^1$-smoothness of $f$ up to the boundary. This equires
an application of the estimates of the Kobayashi-Royden metric given
in Section~3 and the scaling method due to S.Pinchuk; we adapt this to the
almost complex case.

In Section~3 we reduced the problem to the following local
situation. Let $D$ and $D'$ be domains in $\C^2$, $\Gamma$ and
$\Gamma'$ be open $\CC^{\infty}$-smooth pieces of their boundaries,
containing the origin. We assume that an almost complex structure $J$
is defined and $\CC^{\infty}$-smooth in a neighborhood of the closure
$\bar D$, $J(0) = J_{st}$ and $J$ has a diagonal form in a
neighborhood of the origin: $J(z) = diag(a_{11}(z),a_{22}(z))$.
Similarly, we assume that $J'$
is diagonal in a neighborhood of the origin, $J'(z) =
diag(a_{11}'(z),a_{22}'(z))$ and $J'(0) = J_{st}$. The hypersurface
$\Gamma$ (resp. $\Gamma'$) is supposed to be strictly $J$-pseudoconvex
(resp. strictly $J'$-pseudoconvex). Finally, we assume that $f: D
\rightarrow D'$ is a $(J,J')$-biholomorphic map, $1/2$-H\"older
homeomorphism between  $D \cup \Gamma$ and $D' \cup \Gamma'$, such that
$f(\Gamma) = \Gamma'$ and $f(0) = 0$. Finally according to Section~2,
$\Gamma$ is defined in a neighborhood of the origin
by the equation $\rho(z) = 0$ where $\rho(z) = 2Re
z^2 + 2Re K(z) + H(z) + o(\vert z \vert^2)$ and $K(z) = \sum
K_{\mu\nu}z^{\mu\nu}$, $H(z) = \sum h_{\mu\nu}z^{\mu}\bar
z^{\nu}$, $k_{\mu\nu} = k_{\nu\mu}$, $h_{\mu\nu} = \bar
h_{\nu\mu}$. The crucial point is that $H(z^1,0)$ is a positive
hermitian form on $\C$, meaning that in these coordinates $\Gamma$ is
strictly pseudoconvex at the origin with respect to the standard
structure of $\C^2$ (see Lemma~\ref{PP} for the proof). Of course,
$\Gamma'$ admits a similar local representation. In what follows we
assume that we are in this setting. 

Let $(p^k)$ be a sequence of points in $D$  converging to $0$ and let
$\Sigma := \{ z \in \C^2: 2Re z^2 + 2Re K(z^1,0) + H(z^1,0) < 0\}$,
$\Sigma' := \{ z \in \C^2: 2Re z^2 + 2Re K'(z^1,0) + H'(z^1,0) < 0\}$.
The scaling procedure associates with the pair $(f,(p^k)_k)$
a biholomorphism $\phi$ (with respect to the standard structure $J_{st}$)
between $\Sigma$ and $\Sigma'$. Since $\phi$ is obtained as a limit of a
sequence of biholomorphic maps conjugated with $f$, some of their properties
are related and this can be used to study boundary properties of
$f$ and to prove that its cotangent lift is continuous up to the conormal
bundle $\Sigma(\partial D)$.

\subsection{Fixing suitable local coordinates and dilations.}
For any boundary point $t \in
\partial D$ we consider the change of variables $\alpha^t$ defined by 
$$
(z^1)^* = \frac{\partial \rho}{\partial \bar z^2}(t)(z^1 - t^1)
- \frac{\partial \rho}{\partial \bar z^1}(t)(z^2 - t^2),\
(z^2)^* = \sum_{j=1}^2 \frac{\partial \rho}{\partial z^j}(t)(z^j -
t^j).
$$
Then $\alpha^t$ maps $t$ to $0$. The real
normal at $0$ to $\Gamma$ is mapped by $\alpha^t$ to the line $\{
z^1 = 0, y_2 = 0 \}$. 
For every $k$, we denote by $t^k$ the projection of
$p^k$ onto $\partial D$ and by $\alpha^k$ the change of variables
$\alpha^t$ with $ t = t^k$.  Set $\delta_k = dist(p^k,\Gamma)$.
Then $\alpha^k(p^k) = (0,-\delta_k)$ and $\alpha^k(D)=\{2Re z^2 + O(\vert z
\vert^2) < 0\}$ near the origin. Since the sequence $(\alpha^k)_k$ converges
to the identity map, the sequence $(\alpha^k)_*(J)$ of almost
complex structures tends to $J$ as $k \rightarrow \infty$. Moreover
there is a sequence $(L^k)$ of linear automorphisms of $\R^4$
such that $(L^k \circ \alpha^k)_*(J)(0) = J_{st}$.
Then $(L^k \circ \alpha^k)(p^k) = (o(\delta_k),-\delta_k')$ with
$\delta_k' \sim \delta_k$ and
$(L^k \circ \alpha^k)(D) = \{Re (z^2 + \tau_k z^1) + O(\vert z \vert^2) < 0\}$
near the origin, with $\tau_k = o(1)$.
Hence there is sequence $(M^k)$ of
$\C$-linear transformations of $\C^2$, converging to the identity, 
such that $(T^k: = M^k \circ L^k
\circ \alpha^k)$ is a sequence of linear transformations converging to the
identity, and $D^k:= T^k(D)$ is defined near the origin by
$D^k=\{\rho_k(z) = Re z^2 + O(\vert z \vert^2) < 0\}$.
Finally $\tilde p_k = T^k(p^k)= (o(\delta_k),\delta_k'' + io(\delta_k))$ with
$\delta_k''\sim \delta_k$. 
We also denote by $\Gamma^k = \{\rho_k = 0 \}$ the
image of $\Gamma$ under $T^k$.
Furthermore, the sequence of almost complex structures
$(J_k:= (T^k)_*(J))$ converges to $J$ as $k \rightarrow \infty$
and $J_k(0) = J_{st}$.

We proceed quite similarly for the target domain $D'$. 
For $s \in \Gamma'$ we define the transformation $\beta^s$ by
$$
(z^1)^* = \frac{\partial \rho'}{\partial \bar z^2}(s)(z^1 - s^1)
- \frac{\partial \rho'}{\partial \bar z^1}(s)(z^2 - s^2),
(z^2)^* = \sum_{j=1}^2 \frac{\partial \rho'}{\partial z^j}(s)(z^j -
s^j).
$$

Let $s^k$ be the projection of $q^k:=f(p^k)$ onto $\Gamma'$ and
let $\beta^k$ be the corresponding map $\beta^s$ with $s = s^k$. 
The sequence $(q^k)$ converges
to  $0 = f(0)$ so $\beta^k$ tends to the identity. Considering linear
transformations $(L')^k$ and $(M')^k$, we obtain a sequence $(T'^k)$ of
linear transformations converging to the identity and satisfying the
following properties. The domain 
$(D^k)':= T'^k(D')$ is defined near the origin by
$(D^k)'=\{\rho_k'(z) := Re z^2 + O(\vert z \vert^2) < 0\}$,
$\Gamma_k' = \{ \rho_k' = 0 \}$ and $\tilde q_k = T'^k(q^k) =
(o(\varepsilon_k),\varepsilon_k''+ io(\varepsilon_k))$
with $\varepsilon_k'' \sim \varepsilon_k$, where
$\varepsilon_k = dist(q^k,\Gamma')$. 
The sequence of almost complex structures $(J_k':= (T'^k)_*(J'))$
converges to $J'$ as $k \rightarrow \infty$ and $J_k'(0) = J_{st}$.  

Finally, the map $ f^k:= T'^k \circ f \circ (T^k)^{-1}$ satisfies
$f^k(\tilde p_k) = \tilde q_k$  and is
a biholomorphism between the domains $D^k$ and $(D')^k$ 
with respect to the almost
complex structures $J_k$ and $J_k'$. 

Consider now the non isotropic dilations
$\phi_k : (z^1,z^2) \mapsto (\delta_k^{1/2}z^1,\delta_kz^2)$ and
$\psi_k(z^1,z^2)=
(\varepsilon_k^{1/2}z^1,\varepsilon_kz^2)$ and set $\hat f^k =
(\psi_k)^{-1} \circ f^k \circ \phi_k$.
Then the map $\hat f^k$ is biholomorphic with respect to the almost complex
structures $\hat J_k:=((\phi_k)^{-1})_*(J_k)$ and
$\hat J'_k:= (\psi_k^{-1})_*(J'_k)$.
Moreover if $\hat D^k:=\phi_k^{-1}(D^k)$ and
$(\hat{D'})^k:=\psi_k^{-1}((D')^k)$ then
$\hat D^k = \{ z \in \phi_k^{-1}(U): \hat \rho_k(z) < 0\}$
where
$$
\hat \rho_k(z): = \delta_k^{-1}\rho(\phi_k(z)) = 2Re z^2 + \delta_k^{-1}[2
Re K(\delta_k^{1/2}z^1,\delta_kz^2) + H(\delta_k^{1/2}z^1,\delta_kz^2)
+  o(\vert (\delta_k^{1/2}z^1,\delta_kz^2) \vert^2).
$$
and $(\hat D')^k=\{ z \in \phi_k^{-1}(U): \hat \rho'_k(z) < 0\}$
where
$$
\hat \rho'_k(z): = \varepsilon_k^{-1}\rho'(\psi_k(z)) = 2Re z^2 +
\varepsilon_k^{-1}[2 Re K'(\varepsilon_k^{1/2}z^1,\varepsilon_kz^2) +
H'(\varepsilon_k^{1/2}z^1,\varepsilon_kz^2)
+  o(\vert (\varepsilon_k^{1/2}z^1,\varepsilon_kz^2) \vert^2).
$$
Since $U$
is a fixed neighborhood of the origin, the pull-backs $\phi_k^{-1}(U)$
tend to $\C^2$ and the functions $\hat\rho_k$ tend
to $\hat \rho(z) = 2Re z^2 + 2Re K(z^1,0) + H(z^1,0)$ in the $\CC^2$ norm
on any compact subset of $\C^2$. Similarly, since $U'$
is a fixed neighborhood of the origin, the pull-backs $\psi_k^{-1}(U')$
tend to $\C^2$ and the functions $\hat\rho_k'$ tend
to $\hat \rho'(z) = 2Re z^2 + 2Re K'(z^1,0) + H'(z^1,0)$ in the $\CC^2$ norm
on any compact subset of $\C^2$. If $\Sigma :=
\{ z \in \C^2: \hat \rho(z) < 0 \}$ and $\Sigma' := \{ z \in \C^2:
\hat \rho'(z) < 0 \}$ then the sequence of points $\hat p^k =
\phi_k^{-1}(\tilde p_k) \in \hat D^k$ converges to the point $(0,-1) \in
\Sigma$ and the sequence of points $\hat q^k =
\psi^{-1}_k(\tilde q^k) \in \hat{D'}^k$ converges to $(0,-1) \in
\Sigma'$. Finally $\hat{f}^k(\hat p^k) = \hat q^k$.

\subsection{Convergence of the dilated families.} We begin with the following
 
\begin{lemma}\label{convseq}
The sequences $(\hat J'_k)$ and $(\hat J_k)$ of almost complex structures
converge to the standard structure uniformly (with all partial
derivatives of any order) on compact subsets of $\C^2$.
\end{lemma}

\noindent{\it Proof of Lemma~\ref{convseq}.}
Denote by $a_{\nu\mu}^k(z)$ the elements of the matrix
$J_k$. Since $J_k \rightarrow J$ and $J$ is diagonal, we have $a_{\nu\mu}^k
\rightarrow a_{\nu\mu}$ for $\nu = \mu$ and $a_{\nu\mu}^k
\rightarrow 0$ for $\nu \neq \mu$. Moreover, since $J_k(0) =
J_{st}$, $a_{\nu\mu}^k(0) = i$ for $\nu = \mu$ and $a_{\nu\mu}(0) = 0$
for $\nu \neq \mu$.
 The elements $\hat
a_{\nu\mu}^k$ of the matrix
$\hat J_k$ are given by: $\hat a_{\nu\mu}^k(z^1,z^2) = a_{\nu
\mu}^k(\delta_k^{1/2}z^1,\delta_k z^2)$ for $\nu = \mu$, $\hat
a_{12}^k(z^1,z^2) = \delta_k^{1/2}a(\delta_k^{1/2}z^1,\delta_k z^2)$ and
$\hat a_{21}^k(z^1,z^2) = \delta_k^{-1/2}a_{21}^k(\delta_k^{1/2}z^1,\delta_k
z^2)$. This implies the desired result. \qed

\vskip 0,1cm
The next statement is crucial.
\begin{proposition}
\label{scaling}
The sequence $(\hat f^k )$ (together with all derivatives) is a
relatively compact family (with
respect to the compact open topology) on
$\Sigma$; every cluster point $\hat f$ is
a biholomorphism (with respect to $J_{st}$) between $\Sigma$ 
and $\Sigma'$, satisfying $\hat f(0,-1) = (0,-1)$ and
$(\partial \hat f^2/\partial z^2)(0,-1) = 1$. 
\end{proposition}
\noindent{\it Proof of Proposition~\ref{scaling}.}
{\it Step 1: convergence.} Our proof is based on the method
developped by F.Berteloot-G.Coeur\'e and F.Berteloot \cite{BerCo,Ber}.
Consider a domain
$G \subset \C^2$ of the form $G = \{ z \in W: \lambda(z) = 2Re z^2 +
2Re K(z) + H(z) + o(\vert z \vert^2) < 0 \}$ where $W$ is a
neighborhood of the origin. We assume that an almost complex
structure $J$ is diagonal on $W$ and that the hypersurface 
$\{ \lambda = 0 \}$ is strictly $J$-pseudoconvex at any point.
Given $a \in \C^2$ and $\delta > 0$ 
denote by $Q(a,\delta)$ the non-isotropic ball
$Q(a,\delta ) = \{ z: \vert z^1 - a_1 \vert < \delta^{1/2}, \vert z^2
- a_2 \vert < \delta \}$. Denote also by $d_{\delta}$ the non-isotropic
dilation $d_{\delta}(z^1,z^2) = (\delta^{-1/2}z^1,\delta^{-1}z^2)$. 
\begin{lemma}\label{Control}
There exist positive constants $\delta_0, C, r$ satisfying the
following property : for
every $\delta \leq \delta_0$ and for every $J$-holomorphic disc $g:\Delta
\rightarrow G$ such that $g(0) \in Q(0,\delta)$ we have the
inclusion $g(r\Delta) \subset Q(0,C\delta)$.
\end{lemma}
\noindent{\it Proof of Lemma~\ref{Control}.}
Assume by contradiction that there exist positive sequences $\delta_k
\rightarrow 0$, $C_k \rightarrow +\infty$, a sequence $\zeta_k \in
\Delta$, $\zeta_k \rightarrow 0$ and a sequence $g_k: \Delta
\rightarrow G$ of
$J$-holomorphic discs such that $g_k(0) \in Q(0,\delta_k)$ and
$g_k(\zeta_k) \not\in Q(0,C_k\delta_k)$. Denote by $d_k$ the
dilations $d_{\delta}$ with $\delta = \delta_k$ and consider the
composition $h_k = d_k \circ g_k$ defined on $\Delta$. The
dilated domains $G_k:= d_k(G)$ are defined by $\{ z \in d_k(W):
\lambda_k(z):= \delta_k^{-1}\lambda \circ d_k^{-1}(z) < 0 \}$ and the
sequence $(\lambda_k)$ converges uniformly on compact subsets of
$\C^2$ to $\hat \lambda : z \mapsto 2Re z^2 + 2Re K(z) + H(z^1,0)$. Since $J$
is diagonal, the sequence  of structures $J_k:=(d_k)_*(J)$ converges to
$J_{st}$ in the $\CC^2$ norm on compact subsets of $\C^2$.  

The discs $h_k$ are $J_k$-holomorphic and the sequence $(h_k(0))$ is
contained in $Q(0,1)$; passing to a subsequence we may assume that this
converges to a point $p \in \overline{Q(0,1)}$. On the other hand, the
function $\hat \lambda + A\hat \lambda^2$ is
strictly $J_{st}$-plurisubharmonic on $Q(0,5)$ for a suitable constant
$A > 0$. Since the structures $J_k$ tend
to $J_{st}$, the functions $\lambda_k + A\lambda_k^2$
are strictly $J_k$-plurisubharmonic on $Q(0,4)$ for every $k$ large
enough and their Levi forms admit a uniform
lower bound with respect to $k$.
By Proposition~\ref{lowest0,5} the Kobayashi-Royden infinitesimal
pseudometric on $G_k$ admits the following lower bound~:
$K_{G_k}(z,v) \geq C \vert v \vert$ for any $z \in G_k \cap Q(0,3)$,
$v \in \C^2$,
with a positive constant $C$ independent of $k$. Therefore, there exists a
constant $C' > 0$ such that $\vert \vert \vert (dh_k)_\zeta \vert \vert \vert
\leq C'$ for any $\zeta \in (1/2)\Delta$ satisfying
$h_k(\zeta) \in G_k \cap Q(0,3)$.
On the other hand, the sequence $(\vert h_k(\zeta_k) \vert)$ tends to $+
\infty$. Denote by $[0,\zeta_k]$ the segment 
(in $\C$) joining the origin and $\zeta_k$ and let 
$\zeta_k' \in [0,\zeta_k]$ be the point the closest to the origin such
that $h_k([0,\zeta_k']) \subset G_k \cap
\overline{Q(0,2)}$ and $h_k(\zeta_k') \in \partial Q(0,2)$. Since $h_k(0)
\in Q(0,1)$, we have  $\vert h_k(0) - h_k(\zeta_k') \vert \geq C''$
for some constant $C'' > 0$. Let $\zeta_k' = r_k e^{i\theta_k}$, $r_k \in
]0,1[$. Then 
$$
\vert h_k(0) - h_k(\zeta_k') \vert \leq \int_{0}^{r_k}
\vert \vert \vert (dh_k)_{te^{i\theta_k}} \vert \vert \vert dt \leq C'r_k
\rightarrow 0.
$$
This contradiction proves Lemma~\ref{Control}. \qed

\vskip 0,1cm
The statement of Lemma~\ref{Control} remains true if we replace the unit
disc $\Delta$ by the unit ball $\B_2$ in $\C^2$ equipped with an almost
complex structure $\tilde J$ close enough (in the $\CC^2$ norm) to
$J_{st}$. For the proof it is sufficient to foliate $\B_2$ by $\tilde
J$-holomorphic curves through the origin (in view of a smooth
dependence on  small perturbations of $J_{st}$ 
such a foliation is a small
perturbation of the foliation by complex lines through the origin,
see \cite{NiWo}) and apply Lemma~\ref{Control} to the foliation. 

\vskip 0,1cm
As a corollary we have the following
\begin{lemma}
\label{conv}
Let $(M,\tilde J)$ be an almost complex manifold and let $F^k: M
\rightarrow G$ be a sequence of $(\tilde J,J)$-holomorphic maps.
Assume that for some point $p^0 \in M$ we have $F^k(p) = (0,-\delta_k)$,
$\delta_k \rightarrow 0$, and that the sequence $(F^k)$ converges
to $0$ uniformly on compact subsets of $M$. 
Consider the rescaled maps $d_k \circ
F^k$. Then for any compact subset $K \subset M$ the sequence of norms
$(\| d_k \circ F^k \|_{\mathcal C^0(K)})$ is bounded.
\end{lemma}
\noindent{\sl Proof of Lemma~\ref{conv}}. It is sufficient to consider
a covering of a compact subset of $M$ by sufficiently small balls,
similarly to \cite{BerCo}, p.84.
Indeed, consider a covering of $K$ by the balls $p^j + r\B$, $j=0,...,N$
where $r$ is given by Lemma~\ref{Control} and $p^{j+ 1} \in p^j + r\B$ for
any $j$. For $k$ large enough, we obtain
that $F^k(p^0 + r\B) \subset Q(0,2C\delta_k)$, and $F^k(p^1 + r\B)
\subset Q(0,4C^2\delta_k)$. Continuing this process we obtain that
$F^k(p^N + r\mathbb B) \subset Q(0,2^NC^N\delta_k)$.
This proves Lemma~\ref{conv}. \qed

\vskip 0,1cm
Now we return to the proof of Proposition~\ref{scaling}. Lemma
\ref{conv} implies that the sequence $(\hat f^k)$ is bounded (in the $\CC^0$
norm) on any
compact subset $K$ of $\Sigma$. Covering $K$ by small bidiscs,
consider two transversal foliations by $J$-holomorphic curves on every
bidisc. Since the restriction of $\hat f^k$ on every such curve is
uniformly bounded in the $\CC^0$-norm, 
it follows by the elliptic
estimates that this is bounded in $\CC^l$ norm for every $l$ (see
\cite{Si}). Since the bounds are uniform with respect to curves,
the sequence $(\hat f^k)$ is bounded in every
$\CC^l$-norm. So the family $(\hat f^k)$ is relatively compact. 

{\it Step 2: Holomorphy of the limit maps.} Let $(\hat f^{k_s})$ be a
subsequence converging to a smooth map $\hat f$. 
Since $f^{k_s}$ satisfies the holomorphy condition
$\hat J'_{k_s} \circ d\hat f^{k_s} = d\hat f^{k_s} \circ
J_{k_s}$, since $\hat J_{k_s}$ and $\hat J'_{k_s}$ converge to $J_{st}$,
we obtain, passing to the limit in the holomorphy condition, that $\hat f$
is holomorphic with respect to $J_{st}$.

{\it Step 3: Biholomorphy of $\hat f$.}
Since $\hat f(0,-1) = (0,-1) \in \Sigma'$ and
$\Sigma'$ is defined by a plurisubharmonic function, it follows by the
maximum principle that $\hat f(\Sigma) \subset \Sigma'$ (and not just a
subset of $\bar\Sigma'$). Applying a similar argument to the
sequence $(\hat f^k)^{-1}$ of inverse map, we obtain that this converges
(after extraction of a subsequence) to the inverse of $\hat f$. 

Finally the domain $\Sigma$ (resp. $\Sigma'$) is
biholomorphic to ${\mathbb H}$ by means of the transformation $(z^1,z^2)
\mapsto (z^1,z^2 + K(z^1,0))$ (resp. $(z^1,z^2) \mapsto (z^1,z^2 +
K'(z^1,0))$). Since a biholomorphism of ${\mathbb H}$ fixing the point
$(0,-1)$ has the form
$(e^{i\theta}z^1,z^2)$ (see, for instance, \cite{Co}), $\hat f$ is conjugated
to this transformation by the above quadratic biholomorphisms of
$\C^2$. Hence~: 

\begin{eqnarray}
\label{derivative}
\frac{\partial \hat f^2}{\partial z^2}(0,-1) = 1.
\end{eqnarray}

\vskip 0,1cm
\noindent This property will be used in the next Section. \qed

\section{Boundary behavior of  the tangent map} 

We suppose that we are in the local situation described at the
beginning of the previous section. Here we prove two statements
concerning the boundary behavior of the tangent map of $f$ near
$\Gamma$. They are obvious if $f$ is of class $\CC^{1}$ up to
$\Gamma$. In the general situation, their proofs
require the scaling method of the previous section.
Let $p \in \Gamma$. After a local change of coordinates $z$ we
may assume that $p = 0$, $J(0) = J_{st}$ and $J$ is assumed to be diagonal.
In the $z$ coordinates, we 
consider a base $X$ of (1,0) (with respect to $J$) vector fields defined in
Subsection~3.3. Recall that $X_2 = \partial
/\partial z^2 + a(z) \partial/\bar\partial z^2$, $a(0) = 0$,
$X_1(0) = \partial/\partial z^1$ and at every point $z^0$, $X_1(z^0)$
generates the holomorphic tangent space $H_z^J(\partial D - t)$, $t \geq 0$.
If we return to the initial coordinates and move the point $p \in \Gamma$,
we obtain for every $p$ a basis $X_p$ of $(1,0)$ vector fields, defined in a
neighborhood of $p$. Similarly, we define the basis $X'_q$ for
$q \in \partial D'$. 

The elements of the matrix of the tangent map
$df_z$ in the bases $X_p(z)$ and $X'_{f(p)}(z)$ are denoted by
$A_{js}(p,z)$. According to Proposition~\ref{matrix} the function
$A_{22}(p,\cdot)$ is upper bounded on $D$. 

\begin{proposition}
\label{REALITY}
We have:
\begin{itemize}
\item[(a)] Every cluster point of the function $z \mapsto A_{22}(p,z)$
(in the notation of Proposition~\ref{matrix}) is real when $z$ tends to a
point $p \in \partial D$.
\item[(b)] For $z \in D$, let $p \in \Gamma$ such that
$|z-p| = dist(z,\Gamma)$. There exists a constant $A$, independent of
$z \in D$, such that $\vert A_{22}(p,z) \vert \geq A$.
\end{itemize}
\end{proposition}

The proof of these statements use the above scaling
construction. So we use the notations of the previous section. 

\vskip 0,1cm
\noindent{\it Proof of Proposition~\ref{REALITY}}.
(a) Suppose that there exists a sequence of points $(p^k)$ converging
to a boundary point $p$ such that $A_{22}(p,\cdot)$ tends to a complex number
$a$. Applying the above scaling construction,
we obtain a sequence of maps $(\hat f^k)_k$.
Consider the two basis $\hat X^k:= 
\delta_k^{1/2}((\phi_k^{-1}) \circ T^k)(X_1),
\delta_k((\phi_k^{-1})\circ T^k)(X_2))$ and $(\hat
X')^k:= (\varepsilon_k^{-1/2}((\psi_k^{-1}) \circ T'^k)(X'_1),
\varepsilon_k^{-1}((\psi_k^{-1})\circ T'^k)(X'_2))$. These vector
fields tend to the standard (1,0) vector field base of $\C^2$ as $k$
tends to $\infty$. Denote by $\hat A^k_{js}$ the elements of the
matrix of $d\hat f^k(0,-1)$. Then $A^k_{22} \rightarrow (\partial
\hat f^2/\partial z^2)(0,-1) = 1$, according to (\ref{derivative}). On the
other hand, $A^k_{22} = \varepsilon_k^{-1}\delta_k A_{22}$ and tends
to $a$ by the boundary distance preserving property (Proposition~\ref{equiv}).
This gives the statement.

(b) Suppose that there is a sequence of points $(p^k)$ converging
to the boundary such that $A_{22}$ tends to $0$. Repeating precisely
the argument of (a), we obtain that $(\partial \hat f^2/\partial
z^2)(0,-1) = 0$; this contradicts~(\ref{derivative}). \qed

\vskip 0,1cm
In order to establish the next proposition, it is convenient to associate
a wedge with the totally real part of the conormal bundle
$\Sigma_J(\partial D)$ of $\partial D$ as edge.
Consider in $\R^{4} \times \R^{4}$ the set $S = \{ (z,L):
dist((z,L),\Sigma_J(\partial D)) \leq dist(z,\partial D), z \in D \}$.
Then, in a neighborhood $U$ of any totally real point of
$\Sigma_J(\partial D)$, the set S contains a wedge $W_U$ with
$\Sigma_J(\partial D) \cap U$ as totally real edge.

\begin{proposition}
\label{cluster}
Let $K$ be a compact subset of the totally real part of the conormal
bundle $\Sigma_J(\partial D)$. Then the cluster set of the cotangent lift
$\tilde f $ of $f$ on the conormal bundle 
$\Sigma(\partial D)$, when $(z,L)$ tends to $\Sigma_J(\partial D)$
along the wedge $W_U$, is relatively compactly contained 
in the totally real part of $\Sigma(\partial D')$.
\end{proposition}
\noindent{\sl Proof of Proposition~\ref{cluster}}.
Let $(z^k,L^k)$ be a sequence in $W_U$
converging to $(0,\partial_J\rho(0)) = (0,dz^2)$. Set $g = f^{-1}$.
We shall prove that the
sequence of linear forms $Q^k :=
{}^tdg(w^k)L^k$, where $w^k = f(z^k)$, converges to a linear form which up to
a {\it real} factor (in view of Part (a) of
Proposition \ref{REALITY})
coincides with $\partial_J \rho'(0)= dz^2$
(we recall that ${}^t$ denotes the transposed map).
It is sufficient to prove that the first component of $Q^k$ with
respect to the dual basis $(\omega_1,\omega_2)$ of $X$ tends to $0$
and 
the second one is
bounded below from the origin as $k$ tend to infinity. The map $X$ being
of class $\CC^1$ we can replace $X(0)$ by $X(w^k)$.
Since $(z^k,L^k) \in W_U$, we have $L^k
= \omega_2(z^k) + O(\delta_k)$, where $\delta_k$ is the distance from
$z^k$ to 
the boundary. Since $\vert\vert\vert dg_{w^k} \vert\vert\vert =
0(\delta_k^{-1/2})$, we
have $Q^k = {}^tdg_{w^k}(\omega_2(z^k)) + O(\delta_k^{1/2})$.
By Proposition~\ref{matrix}, the
components of ${}^tdg_{w^k}(\omega_2(z^k))$ with respect to the basis
$(\omega_1(z^k),\omega_2(z^k))$ are the elements of the 
second line of the matrix 
$dg_{w^k}$ with respect to the basis $X'(w^k)$ and $X(z^k)$. So its first 
component is $0(\delta_k^{1/2})$ and tends to $0$ as $k$ tends to
infinity. Finally the component $A_{22}^k$ is bounded below from the
origin by Part (b) of Proposition~\ref{REALITY}. \qed

\vskip 0,1cm
\noindent{\it Proof of Theorem \ref{MainTheorem}}. In view of
Proposition~\ref{cluster}, we may apply Proposition~\ref{wedges} to the
cotangent lift $\tilde f$ of $f$.
This gives the statement of Theorem~\ref{MainTheorem}. \qed

\section{Appendices}

\subsection{Appendix 1 : Uniform estimates of the Kobayashi-Royden metric}

We prove Proposition~\ref{lowest0,5}, restated in Proposition~\ref{addest}.
Our method is based on Sibony's approach \cite{Si}.

For our construction we need plurisubharmonic functions with
logarithmic singularities on almost complex manifolds. The function
$\log \vert z \vert$ is not $J$-plurisubharmonic in a neighborhood of
the origin even if the structure $J$ is $\CC^2$ close to the standard
one. However, for a suitable positive constant $A > 0$ the function
$\log \vert z \vert + A\vert z \vert$ is $J$-plurisubharmonic on the
unit ball $\B$ for any almost complex structure $J$ with $\| J
- J_{st}\|_{C^{2}(\bar\B)}$ small enough. This useful
observation due to E.Chirka can be easily established by direct
computation of the Levi form, see \cite{GaSu} (we point out that the
Levi form of $\vert z \vert$ goes to $+\infty$ at the origin
neutralizing the growth of the logarithm). This implies the following~: 

\begin{lemma}\label{lemlem}
Let $r < 1$ and let $\theta_r$ be a 
smooth nondecreasing function on
$\R^+$ such that $\theta_r(s)= s$ for $s \leq r/3$ and $\theta_r(s) =
1$ for $s \geq 2r/3$. Let $(M,J)$ be an almost complex manifold, and
let $p$ be a point of $M$. Then there exists a neighborhood $U$ of
$p$, positive constants $A = A(r)$, $B=B(r)$ and a diffeomorphism $z:U
\rightarrow  B$ such that $z(p) = 0$, $dz(p) \circ J(p) \circ
dz^{-1}(0) = J_{st}$ and the function ${\rm log}(\theta_r(\vert z
\vert^2)) + \theta_r(A\vert z \vert) + B\vert z \vert^2$ is
$J$-plurisubharmonic on $U$.
\end{lemma}

The main estimate of the Kobayashi-Royden metric is given by the following

\begin{proposition}\label{addest}
Let $D$ be a domain in an almost complex
manifold $(M,J)$, let $p \in \bar{D}$, let $U$ be a neighborhood of
$p$ in $M$ (not necessarily contained in $D$) and let $z:U \rightarrow
\B$ be the diffeomorphism given by Lemma~\ref{lemlem}.  
Let $u$ be a $\mathcal C^2$ function on $\bar{D}$, negative and
$J$-plurisubharmonic on $D$. We assume that $-L \leq u < 0$ on $D \cap
U$ and that $u-c|z|^2$ is $J$-plurisubharmonic on $D \cap U$, where
$c$ and $L$ are positive constants. Then there exists a neighborhood
$U'$ of $p$ and a constant $c'
> 0$, depending on $c$ and $L$ only, such that :

\begin{equation}\label{e1}
K_{(D,J)}(q,v) \geq c'\frac{\|v\|}{|u(q)|^{1/2}},
\end{equation}
for every $q \in D \cap U'$ and every $v \in T_qM$.
\end{proposition}

\noindent{\it Proof of proposition~\ref{addest}.}
{\it Step 1: Local hyperbolicity.} We prove the following
rough estimate 
\begin{equation}\label{e2}
K_{(D,J)}(q,v) \geq s K_{(D \cap U,J)}(q,v)
\end{equation}
which allows to localize the proof ($s$ is a positive constant).
Let $0<r<1$ be such that the set $V_1:=\{q \in U :
|z(q)| \leq \sqrt{r}\}$ is relatively compact in $U$ and let $\theta_r$ 
be a smooth nondecreasing function on $\R^+$ such that $\theta_r(s)= s$ for
$s \leq r/3$ and $\theta_r(s) = 1$ for $s \geq 2r/3$. According to Lemma
~\ref{lemlem}, there exist uniform positive constants $A$ and $B$ such
that the function 
$
{\rm log}(\theta_r(|z-z(q)|^2))+ \theta_r(A|z-z(q)|)+ B|z|^2
$
is $J$-plurisubharmonic on $U$ for every $q \in V$. 
By assumption the function
$u-c|z|^2$ is $J$-plurisubharmonic on $D \cap U$. Set $\tau=2B/c$ and
define, for every point $q \in V$, the function $\Psi_{q}$ by~:
$$
\left\{
\begin{array}{lll}
\Psi_{q}(z) &=& \theta_r(|z-z(q)|^2)\exp(\theta_r(A|z-z(q)|)) 
\exp(\tau u(z))\ {\rm if} \  z \in D \cap U,\\
& & \\
\Psi_{q} &=& \exp(1+\tau u) \ {\rm on} \ D \backslash U.
\end{array}
\right.
$$

Then for every $0 < \varepsilon \leq B$, the function ${\rm
log}(\Psi_{q})-\varepsilon|z|^2$ is $J$-plurisubharmonic on $D \cap U$
and hence $\Psi_{q}$ is $J$-plurisubharmonic on $D \cap U$. Since
$\Psi_{q}$ coincides with $\exp(\tau u)$ outside $U$, it is globally
$J$-plurisubharmonic on $D$. 

Let $f \in \mathcal O_{J}(\Delta,D)$ be such that $f(0)=q \in V_1$ and
$(\partial f/\partial x)(0) = v/\alpha$ where $v \in T_qM$ and $\alpha
>0$. For $\zeta$ sufficiently close to 0 we have
$f(\zeta) = q + df_0(\zeta) + \mathcal O(|\zeta|^2)$.
Setting $\zeta= \zeta_1+i\zeta_2$ and using
the $J$-holomorphy condition $df_0\circ J_{st} = J \circ
df_0$, we may write $
df_0(\zeta) = \zeta_1 df_0(\partial /
\partial x) + \zeta_2 J(df_0(\partial / \partial x))
$.
Consider the function
$
\varphi(\zeta) = \Psi_q(f(\zeta))/|\zeta|^2
$
which is subharmonic on
$\Delta \backslash \{0\}$. Since
$
\varphi(\zeta) = |f(\zeta)-q|^2/|\zeta|^2 \exp(A|f(\zeta)-q|) 
\exp(\tau u(f(\zeta)))
$
for $\zeta$ close to 0 and 
$
\|df_0(\zeta)\| \leq |\zeta| (\|I+J\|\,\|df_0(\partial
/\partial x)\|)
$,
we obtain that $\limsup_{\zeta \rightarrow 0}\varphi(\zeta)
$ is finite. Moreover setting $\zeta_2=0$ we have 
$
\limsup_{\zeta \rightarrow 0}\varphi(\zeta) \geq \|df_0(\partial
/\partial x)\|^2\exp(-2B|u(q)|/c).
$
Applying the maximum principle to a subharmonic extension of $\varphi$
on $\Delta$ we obtain the inequality 
$
\|df_0(\partial / \partial x)\|^2 \leq \exp(1+2B|u(q)|/c).
$
Hence, by definition of the Kobayashi-Royden infinitesimal pseudometric, 
we obtain for every $q \in D \cap V_1$, $v \in T_qM$~:
$$
K_{(D,J)}(q,v) \geq \left(\exp\left(-1-2B\frac{|u(q)|}{c}\right)
\right)^{1/2}\|v\|.
$$
We denote by $d_{(M,J)}^K$ the integrated pseudodistance of the
Kobayashi-Royden infinitesimal pseudometric.
According to the almost complex version of
Royden's theorem \cite{kr99}, it coincides with the usual Kobayashi
pseudodistance on $(M,J)$ defined by means of $J$-holomorphic discs.
Consider now the Kobayashi ball $B_{(D,J)}(q,\alpha)=\{w \in D :
d_{(D,J)}^K(w,q)<\alpha\}$. It follows from Lemma~2.2 of \cite{ChCoSu} (whose
proof is identical in the almost complex setting) that
there is a neighborhood $V$ of $p$, relatively compact in $V_1$ 
and a positive constant $s<1$,
independent of $q$, such that for every $f \in \mathcal
O_{J}(\Delta,D)$ satisfying $f(0) \in D \cap V$ we have $f(s\Delta)
\subset D \cap U$. This gives the inequality (\ref{e2}).

{\it Step 2}. It follows from 
(\ref{e2}) that there is a 
neighborhood $V$ of $p$ in $\C^n$, contained in $U$
and a positive constant $s$ such that 
$D_{(D,J)}(q,v) \geq s K_{(D\cap U,J)}(q,v)$ 
for every $q \in V,\ v \in T_qM$. 
Consider a positive constant $r$ that will be specified later and
let $\theta$ be a smooth nondecreasing function
on $\R^+$ such that $\theta(x) =x$ for $x \leq 1/3$ and $\theta(x) =
1$ for $x \geq 2/3$. Restricting $U$ if necesssary it follows from 
Lemma~\ref{lemlem} that the function 
$
\log(\theta(|(z-q)/r|^2)) + A|z-q| + B|(z-q)/r|^2
$
is $J$-plurisubharmonic on $D \cap U$, independently of $q$ and $r$.

Consider now the function 
$
\Psi_{q}(z) =
\theta\left(\frac{|z-q|^2}{r^2}\right)\exp(A|z-q|)\exp(\tau u(z))
$
where $\tau=1/|u(q)|$ and $r=(2B|u(q)|/c)^{1/2}$. 
Since the function $\tau u - 2B |(z-q)/r|^2$ is $J_{st}$-plurisubharmonic,
we may assume, shrinking $U$ if necessary, that the function
$\tau u - B |(z-q)/r|^2$ is $J$-plurisubharmonic on $D \cap U$. Hence the
function $\log(\Psi_q)$ is $J$-plurisubharmonic on $D \cap U$.
Let $q \in V$, let $v \in T_qM$ and let 
$f:\Delta \rightarrow D$ be a $J$-holomorphic map be such that
$f(0) = q$ and $df_0(\partial / \partial x) = v/\alpha$ where $\alpha
>0$. We have $f(\zeta) = q + df_0(\zeta) + \mathcal
O(|\zeta|^2)$. Setting $\zeta= \zeta_1+i\zeta_2$ and using the
$J$-holomorphy condition $df_0\circ J_{st} = J \circ df_0$, we may
write 
$
df_0(\zeta) = q +\zeta_1 df_0(\partial /
\partial x) + \zeta_2 J(df_0(\partial / \partial x))
$.
Consider the function
$
\varphi(\zeta) = \Psi_q(f(\zeta))/|\zeta|^2
$
which is subharmonic on
$\Delta\backslash \{0\}$.  Since 
$
\varphi(\zeta) = |f(\zeta)-q|^2/(r^2|\zeta|^2) \exp(\tau u(f(\zeta)))
$
and 
$
|df_0(\zeta)| \leq |\zeta| (\|I+J\|\,\|df_0(\partial
/\partial x)\|
$, 
we obtain that $\limsup_{\zeta \rightarrow 0}\varphi(\zeta)$
is finite.
Setting $\zeta_2=0$ we obtain $
\limsup_{\zeta \rightarrow
0}\phi(\zeta) \geq \|v\|^2\exp(2)/(r^2\alpha^2)$.
There exists a positive constant $C'$, independent of $q$, such
that $|z-q| \leq C'$ on $D$.  Applying the maximum
principle to a subharmonic extension of $\phi$ on $\Delta$, we obtain the
inequality
$$
\alpha \geq \sqrt{\frac{c}{2B\exp(1+AC')}}\|v\|^2/|u(q)|^{1/2}.
$$
This completes the proof. \qed

\subsection{Appendix 2 : Canonical lift of an almost complex structure
to the cotangent bundle}
We recall the definition of the canonical lift of an almost
complex structure $J$ on $M$ to the cotangent bundle $T^*M$, following~
\cite{YI}. Set $m=
2n$. We use the following notations. Suffixes A,B,C,D take the values
$1$ to $2m$, suffixes $a,b,c,\dots$,$h,i,j,\dots$ take the values $1$ to
$m$ and $\bar j = j+ m$, $\dots$ The summation notation for
repeated indices is used. If the notation $(\varepsilon_{AB})$,
$(\varepsilon^{AB})$, $(F_{B}^{{}A})$ is used for matrices, the suffix
on the left indicates the column and the suffix on the right indicates
the row. We denote local coordinates on $M$ by $(x^1,\dots,x^n)$ and by
$(p_1,\dots,p_n)$ the fiber coordinates.

Recall that the cotangent space $T^*(M)$ of $M$ possesses the {\it
canonical contact form} $ \theta$ given in local coordinates by
$\theta = p_idx^i.$
The cotangent lift $\varphi^*$ of any diffeomorphism $\varphi$ of $M$
is contact with respect to $\theta$, that is $\theta$ does not depend on
the choice of local coordinates on $T^*(M)$. 

The exterior derivative $d\theta$ of $\theta$ defines {\it the
canonical
symplectic structure} of $T^*(M)$:
$d\theta = dp_i \wedge dx^i$
which is also independent of local coordinates in view of the
invariance of the exterior derivative. Setting $d\theta =
(1/2)\varepsilon_{CB}dx^C \wedge dx^B$ (where $dx^{\bar j} =
dp_j$), we have

$$
(\varepsilon_{CB}) = \left(
\begin{array}{cc}
0 & I_n \\
-I_n  & 0
\end{array}
\right).
$$

Denote by $(\varepsilon^{BA})$ the inverse matrix and write
$\varepsilon^{-1}$ for the tensor field of type (2,0) whose component
are $(\varepsilon^{BA})$. By construction, this definition does not
depend on the choice of local coordinates.

Let now $E$ be a tensor field of type (1,1) on $M$. If $E$ has
components $E_i^{\, h}$ and $E_i^{*h}$ relative to local coordinates $x$
and $x^*$ repectively, then
$p_a^*E_{i}^{*\,a} = p_aE_j^{\, b}\frac{\partial x^j}{\partial x^{*i}}.$
If we interpret a change of coordinates as a diffeomorphism $x^* =
x^*(x) = \varphi(x)$ we denote by $E^*$ the direct image of the tensor $E$
under the action of $\varphi$. In the case where $E$ is an almost
complex structure (that is $E^2 = -Id$), then $\varphi$ is a biholomorphism
between $(M,E)$ and $(M,E^*)$. Any (1,1) tensor field $E$ on $M$
canonically defines a contact form on $E^*M$ via $
\sigma = p_aE_b^{\, a}dx^b.$
Since $(\varphi^*)^*( p_a^*E_b^{{*}\,a}dx^{*b}) = \sigma,$
$\sigma$ does not depend on a choice of local coordinates (here
$\varphi^*$ is the cotangent lift of $\varphi$). Then this canonically
defines the symplectic form

$$
d\sigma = p_a\frac{\partial E_b^{\, a}}{\partial x^c}dx^c \wedge dx^b
+ E_b^{\, a}dp_a \wedge dx^b.
$$
The cotangent lift $\varphi^*$ of a diffeomorphism $\varphi$ is a
symplectomorphism for $d\sigma$. We may write
$d\sigma = (1/2)\tau_{CB}dx^C \wedge dx^B$
where $x^{\bar i} = p_i$; so we have

$$
\tau_{ji} = p_a \left ( \frac{\partial E_i^{\, a}}{\partial x^j} -
\frac{\partial E_j^{\, a}}{\partial x^i} \right ), \tau_{\bar{j} i} =
E_i^{\, j}, \tau_{j\bar{i}} = -E_j^{\, i}, \tau_{\bar{j}\bar{i}} = 0.
$$

We write $\widehat{E}$ for the tensor field of type (1,1) on $T^*(M)$ whose
components $\widehat{E}_B^{{}A}$  are given by 
$\widehat{E}_B^{{}A} = \tau_{BC}\varepsilon^{CA}.$
Thus $
\widehat{E}_i^{\, h} = E_i^{\, h}, \ \widehat{E}_{\bar{i}}^{\, h} = 0$
and
$\widehat{E}_i^{\, \bar{h}} = 
p_a \left( \frac{\partial E_i^{\, a}}{\partial x^j} -
\frac{\partial E_j^{\, a}}{\partial x^i} \right ), 
\widehat{E}_{\bar{i}}^{\,\bar{h}} = E_h^{\, i}.$
In the matrix form we have

$$
\widehat{E} = \left(
\begin{array}{cll}
E_i^{\, h} & & 0 \\
 p_a \left ( \frac{\partial E_i^{\, a}}{\partial x^j} -
\frac{\partial E_j^{\, a}}{\partial x^i} \right ) & & E_h^{\, i}
\end{array}
\right).
$$

By construction, the complete lift $\widehat{E}$ has the following
{\it invariance property}~: if $\varphi$ is a local diffeomorphism of $M$
transforming $E$ to $E'$, then the direct image of $\widehat{E}$ under the
cotangent lift $\psi: = \varphi^*$ is $\widehat{E'}$.
In general, $\widehat{E}$ is not an almost complex
structure, even if $E$ is.  Moreover, one can show \cite{YI} that
$\widehat{J}$ is a complex structure if and only if $J$ is integrable.
One may however construct an almost complex structure on $T^*(M)$ as follows.

Let $S$ be a tensor field of type (1,s) on $M$. We may consider the
tensor field $\gamma S$ of type $(1,s-1)$ on $T^*M$, defined in local
canonical coordinates on $T^*M$ by the expression

$$
\gamma S = p_aS_{i_s...i_2i_1}^{\,\,\,a}dx^{i_s} \otimes
\cdots \otimes dx^{i_2}\otimes \frac{\partial}{\partial p_{i_1}}.
$$

In particular, if $T$ is a tensor field of type (1,2) on $M$, then $\gamma T$
has components

$$
\gamma T = \left(
\begin{array}{cll}
0 & & 0\\
 p_a T_{ji}^{\,\,a} & & 0 
\end{array}
\right)
$$
in the local canonical coordinates on $T^*M$.

Let $F$ be a (1,1) tensor field on $M$. Its Nijenhuis tensor $N$
is the tensor field of type (1,2) on $M$ acting on two vector fields $X$ and
$Y$ by

$$
N(X,Y) = [FX,FY] - F[FX,Y] - F[X,FY] + F^2[X,Y].
$$

By $NF$ we denote the tensor field acting by $(NF)(X,Y) =
N(X,FY)$. The following proposition is proved in \cite{YI} (p.256).

\begin{proposition}
Let $J$ be an almost complex structure on $M$. Then
\begin{equation}\label{EEE}
\tilde J := \widehat J + (1/2)\gamma(NF)
\end{equation}
is an almost complex structure on the cotangent bundle $T^*(M)$.
\end{proposition}

We stress that the definition of the tensor $\tilde{J}$ is independent of
the choice of coordinates on $T^*M$. Therefore if $\phi$ is a biholomorphism
between two almost complex manifolds $(M,J)$ and $(M',J')$, then its
cotangent lift is a biholomorphism between $(T^*(M),\tilde J)$ and
$(T^*(M'), \tilde J')$. Indeed one can view $\phi$ as a change of coordinates
on $M$, $J'$ representing $J$ in the new coordinates. The cotangent lift
$\phi^*$ defines a change of coordinates on $T^*M$ and $\tilde{J}'$
represents $\tilde{J}$ in the new coordinates.
So the assertion $(i)$ of Proposition~\ref{prop-lift} holds. Property
$(ii)$ of Proposition~\ref{prop-lift} is immediate in view of the definition
of $\tilde{J}$ given by~(\ref{EEE}).

\end{document}